\documentclass[11pt,twoside]{article}
\usepackage{amsmath, amsthm, amscd, amsfonts, amssymb, graphicx, color}

\date{}

\begin{document}

\centerline{}

\centerline {\Large{\bf A study on Best Approximation and Banach Algebra }}
\centerline {\Large{\bf in \,$n$-normed linear space}}

\newcommand{\mvec}[1]{\mbox{\bfseries\itshape #1}}
\centerline{}
\centerline{\textbf{Prasenjit Ghosh}}
\centerline{Department of Mathematics,}
\centerline{Barwan N. S. High School (HS), Barwan, }
\centerline{Murshidabad, 742161, West Bengal, India}
\centerline{e-mail: prasenjitpuremath@gmail.com}
\centerline{}
\centerline{\textbf{T. K. Samanta}}
\centerline{Department of Mathematics, Uluberia College,}
\centerline{Uluberia, Howrah, 711315,  West Bengal, India}
\centerline{e-mail: mumpu$_{-}$tapas5@yahoo.co.in}

\newtheorem{Theorem}{\quad Theorem}[section]

\newtheorem{definition}[Theorem]{\quad Definition}

\newtheorem{theorem}[Theorem]{\quad Theorem}

\newtheorem{remark}[Theorem]{\quad Remark}

\newtheorem{corollary}[Theorem]{\quad Corollary}

\newtheorem{note}[Theorem]{\quad Note}

\newtheorem{lemma}[Theorem]{\quad Lemma}

\newtheorem{example}[Theorem]{\quad Example}

\newtheorem{result}[Theorem]{\quad Result}
\newtheorem{conclusion}[Theorem]{\quad Conclusion}

\newtheorem{proposition}[Theorem]{\quad Proposition}

\begin{abstract}
\textbf{\emph{The idea of best approximation in linear \,$n$-normed space is presented and some examples showing various possibilities of best approximations in linear \,$n$-normed space is given.\,Also, we study strictly convex \,$n$-norm and enquire about the uniqueness of best approximations in \,$n$-normed linear space.\,Furthermore, best approximations in \,$n$-Hilbert space is discussed.\,Moreover, the notion of a Banach algebra in \,$n$-Banach space is presented and some examples are discussed.\,A set-theoretic property of invertible and non-invertible elements in a \,$n$-Banach algebra is explained and then topological divisor of zero in \,$n$-Banach algebra is defined.\,Finally, we introduce the notion of a complex homeomorphism in a \,$n$-Banach algebra and derive Gleason, Kahane, Zelazko type theorem with the help of complex \,$b$-homeomorphism in the case of \,$n$-Banach algebra.}}
\end{abstract}


\section{Introduction}

\smallskip\hspace{.6 cm}
We know that any real-valued continuous function in \,$[\,a,\,b\,]$\, can be approximated by a set of polynomials.\,There are also another facts on approximations, for example any open interval in the real line can be approximated by a set of closed intervals.\,A natural setting for the problem of approximation is as follows.

Let \,$X$\, be a normed linear space and \,$G$\, be a subset of \,$X$.\,Let \,$x_{0} \,\in\, X$\, and \,$\delta \,=\, \inf\limits_{g \,\in\, G}\,\left\|\,x_{0} \,-\, g\,\right\|$.\,Then $\delta$\, is the distance of \,$x_{0}$\, from \,$G$.\,The problem for best approximation of  \,$x_{0}$\, out of the element of \,$G$\, is to find an element \,$g_{0} \,\in\, G$\, such that \,$\delta \,=\, \left\|\,x_{0} \,-\, g_{0}\,\right\|$ .\,We see that a best approximation which occurs for \,$g_{0}$, an element of minimum
distance from the given \,$x_{0}$.\,Such a \,$g_{0} \,\in\, G$\, may or may not exist; which raises the problem of existence.\,The problem of uniqueness is of practical interest too, since for given \,$x_{0}$ and $G$\, there may be more than one best approximation.

Nagumo \cite{MN} introduced the notion of Banach algebra in 1936 and thereafter further development of the theory of Banach algebra was given by Gelfand et al. \cite{IMG, IM}.\,In recent times, various Banach algebra's techniques are used to simplify the theories related to matrices, operators, integral equations and dynamical systems etc.

The idea of linear 2-normed space was first introduced by S. Gahler \cite{Gahler} and thereafter the geometric structure of linear 2-normed spaces was developed by Y. J. Cho and R. W. Freese \cite{Freese}.\,The concept of \,$2$-Banach space is briefly discussed in \cite{White}.\,Mohammed and Siddiqui \cite{NMA} introduced the concept of \,$2$-Banach algebra and derived some known results of the usual Banach algebra in \,$2$-Banach algebra.\,For more on \,$2$-normed algebras one can go through the papers \cite{NSS, RU}.\,H.\,Gunawan and Mashadi \cite{Mashadi} developed a generalization of a linear $2$-normed space for \,$n \,\geq\, 2$.\,Some fundamental results of classical normed space with respect to \,$b$-linear functional in linear\;$n$-normed space  have been studied by P. Ghosh and T. K. Samanta \cite{Prasenjit, K, KK}.\,Also they have studied few fixed point theorems in linear \,$n$-normed space \cite{KP}.\,The concept of \,$2$-inner product space was first introduced by Diminnie et al.\,\cite{Diminnie} in 1970's.\;In 1989, A.\,Misiak \cite{Misiak} developed a generalization of a \,$2$-inner product space for \,$n \,\geq\, 2$.

In this paper, first we have studied best approximations theory in linear \,$n$-normed space and observed that finite dimensionality plays an important role for existence of best approximation.\,Strictly convex and uniformly convex in linear \,$n$-normed space are also discussed.\,The notion of a Banach algebra in \,$n$-Banach spaces is being discussed with a few examples.\,We shall establish a set-theoretic property of invertible and non-invertible elements in a \,$n$-Banach algebra and then present topological divisor of zero in a \,$n$-Banach algebra.\,Finally, the idea of a complex homeomorphism in a \,$n$-Banach algebra is being given and few results with respect to complex \,$b$-homeomorphism in \,$n$-Banach algebra are established.

\section{Preliminaries}
\smallskip\hspace{.6 cm}
In this section, we give some necessary definitions and theorems.

\begin{definition}\cite{Mashadi}
Let \,$X$\, be a linear space over the field \,$ \mathbb{K}$, where \,$ \mathbb{K} $\, is the real or complex numbers field with \,$\text{dim}\,X \,\geq\, n$, where \,$n$\, is a positive integer.\;A real valued function \,$\left \|\,\cdot \,,\, \cdots \,,\, \cdot \,\right \| \,:\, X^{\,n} \,\to\, \mathbb{R}$\, is called an n-norm on \,$X$\, if
\begin{itemize}
\item[(N1)]$\left\|\,x_{\,1} \,,\, x_{\,2} \,,\, \cdots \,,\, x_{\,n}\,\right\| \,=\,0$\, if and only if \,$x_{\,1},\, \cdots,\, x_{\,n}$\, are linearly dependent,
\item[(N2)]$\left\|\,x_{\,1} \,,\, x_{\,2} \,,\, \cdots \,,\, x_{\,n}\,\right\|$\; is invariant under permutations of \,$x_{\,1},\, x_{\,2},\, \cdots,\, x_{\,n}$,
\item[(N3)]$\left\|\,\alpha\,x_{\,1} \,,\, x_{\,2} \,,\, \cdots \,,\, x_{\,n}\,\right\| \,=\, |\,\alpha\,|\, \left\|\,x_{\,1} \,,\, x_{\,2} \,,\, \cdots \,,\, x_{\,n}\,\right\|\; \;\;\forall \;\; \alpha \,\in\, \mathbb{K}$,
\item[(N4)]$\left\|\,x \,+\, y \,,\, x_{\,2} \,,\, \cdots \,,\, x_{\,n}\,\right\| \,\leq\, \left\|\,x \,,\, x_{\,2} \,,\, \cdots \,,\, x_{\,n}\,\right\| \,+\,  \left\|\,y \,,\, x_{\,2} \,,\, \cdots \,,\, x_{\,n}\,\right\|$
\end{itemize}
hold for all \,$x,\, y,\, x_{\,1},\, x_{\,2},\, \cdots,\, x_{\,n} \,\in\, X$.\;The pair \,$\left(\,X,\, \left \|\,\cdot,\, \cdots,\, \cdot \,\right \| \,\right)$\; is then called a linear n-normed space. 
\end{definition}

\begin{definition}\cite{Mashadi}
A sequence \,$\{\,x_{\,k}\,\} \,\subseteq\, X$\, is said to converge to \,$x \,\in\, X$\; if 
\[\lim\limits_{k \to \infty}\,\left\|\,x_{\,k} \,-\, x \,,\, e_{\,2} \,,\, \cdots \,,\, e_{\,n} \,\right\| \,=\, 0\]
for every \,$ e_{\,2},\, \cdots,\, e_{\,n} \,\in\, X$\, and it is called a Cauchy sequence if 
\[\lim\limits_{l \,,\, k \to \infty}\,\left\|\,x_{\,l} \,-\, x_{\,k} \,,\, e_{\,2} \,,\, \cdots \,,\, e_{\,n}\,\right\| \,=\, 0\]
for every \,$ e_{\,2},\, \cdots,\, e_{\,n} \,\in\, X$.\;The space \,$X$\, is said to be complete or n-Banach space if every Cauchy sequence in this space is convergent in \,$X$.\,If a sequence \,$\{\,x_{\,k}\,\}$\, is converges to \,$x$\, then it will be denoted by \,$\lim\limits_{k \,\to\, \infty} x_{\,k} \,=\, x$.   
\end{definition}

\begin{definition}\cite{Soenjaya}
We define the following open and closed ball or sphere in \,$X$: 
\[B_{\,\{\,e_{\,2} \,,\, \cdots \,,\, e_{\,n}\,\}}\,(\,a \,,\, \delta\,) \,=\, \left\{\,x \,\in\, X \,:\, \left\|\,x \,-\, a \,,\, e_{\,2} \,,\, \cdots \,,\, e_{\,n}\,\right\| \,<\, \delta \,\right\}\;\text{and}\]
\[B_{\,\{\,e_{\,2} \,,\, \cdots \,,\, e_{\,n}\,\}}\,[\,a \,,\, \delta\,] \,=\, \left\{\,x \,\in\, X \,:\, \left\|\,x \,-\, a \,,\, e_{\,2} \,,\, \cdots \,,\, e_{\,n}\,\right\| \,\leq\, \delta\,\right\},\hspace{.5cm}\]
where \,$a,\, e_{\,2},\, \cdots,\, e_{\,n} \,\in\, X$\, and \,$\delta$\, be a positive number.
\end{definition}

\begin{definition}\cite{Soenjaya}
A subset \,$G$\, of \,$X$\, is said to be open in \,$X$\, if for all \,$a \,\in\, G $, there exist \,$e_{\,2},\, \cdots,\, e_{\,n} \,\in\, X $\, and \, $\delta \,>\, 0 $\; such that \,$B_{\,\{\,e_{\,2} \,,\, \cdots \,,\, e_{\,n}\,\}}\,(\,a \,,\, \delta\,) \,\subseteq\, G$.
\end{definition}

\begin{definition}\cite{Soenjaya}
Let \,$ A \,\subseteq\, X$.\;Then the closure of \,$A$\, is defined as 
\[\overline{A} \,=\, \left\{\, x \,\in\, X \;|\; \,\exists\, \;\{\,x_{\,k}\,\} \,\in\, A \;\;\textit{with}\;  \lim\limits_{k \,\to\, \infty} x_{\,k} \,=\, x \,\right\}.\]
The set \,$ A $\, is said to be closed if $ A \,=\, \overline{A}$. 
\end{definition}

\begin{definition}\cite{Soenjaya}
A point \,$x \,\in\, X$\, is called a boundary point of a subset \,$A$\, of \,$X$\, if in every open ball with center at \,$x$\, there are points of \,$A$\, as well as points of \,$\overline{\,A}$.\,The boundary of \,$A$\, is the set of all boundary points of \,$A$.
\end{definition}

\begin{definition}\cite{Prasenjit}
Let \,$W$\, be a subspace of \,$X$\, and \,$b_{\,2},\, b_{\,3},\, \cdots,\, b_{\,n}$\; be fixed elements in \,$X$\, and \,$\left<\,b_{\,i}\,\right>$\, denote the subspaces of \,$X$\, generated by \,$b_{\,i}$, for \,$i \,=\, 2,\, 3,\, \cdots,\,n $.\;Then a map \,$T \,:\, W \,\times\,\left<\,b_{\,2}\,\right> \,\times\, \cdots \,\times\, \left<\,b_{\,n}\,\right> \,\to\, \mathbb{K}$\; is called a b-linear functional on \,$W \,\times\, \left<\,b_{\,2}\,\right> \,\times\, \cdots \,\times\, \left<\,b_{\,n}\,\right>$, if for every \,$x,\, y \,\in\, W$\, and \,$k \,\in\, \mathbb{K}$, the following conditions hold:
\begin{itemize}
\item[(I)]\hspace{.2cm} $T\,(\,x \,+\, y,\, b_{\,2},\, \cdots,\, b_{\,n}\,) \,=\, T\,(\,x,\, b_{\,2},\, \cdots,\, b_{\,n}\,) \,+\, T\,(\,y,\, b_{\,2},\, \cdots,\, b_{\,n}\,)$
\item[(II)]\hspace{.2cm} $T\,(\,k\,x,\, b_{\,2},\, \cdots,\, b_{\,n}\,) \,=\, k\,T\,(\,x,\, b_{\,2},\, \cdots,\, b_{\,n}\,)$. 
\end{itemize}
A b-linear functional is said to be bounded if there exists a real number \,$M \,>\, 0$\; such that
\[\left|\,T\,(\,x,\, b_{\,2},\, \cdots,\, b_{\,n}\,)\,\right| \,\leq\, M\; \left\|\,x,\, b_{\,2},\, \cdots,\, b_{\,n}\,\right\|\; \;\forall\; x \,\in\, W.\]
The norm of the bounded b-linear functional \,$T$\, is defined by
\[\|\,T\,\| \,=\, \inf\,\left\{\,M \,>\, 0 \;:\; \left|\,T\,(\,x,\, b_{\,2},\, \cdots,\, b_{\,n}\,)\,\right| \,\leq\, M\,\left\|\,x,\, b_{\,2},\, \cdots,\, b_{\,n}\,\right\|\; \;\forall\; x \,\in\, W\,\right\}.\]
\end{definition}

If \,$T$\, be a bounded \,$b$-linear functional on \,$W \,\times\, \left<\,b_{\,2}\,\right> \,\times\, \cdots \,\times\, \left<\,b_{\,n}\,\right>$, norm of \,$T$\, can be expressed by any one of the following equivalent formula:
\begin{itemize}
\item[$(I)$]\hspace{.2cm}$\|\,T\,\| \,=\, \sup\,\left\{\,\left|\,T\,(\,x,\, b_{\,2},\, \cdots,\, b_{\,n}\,)\,\right| \;:\; \left\|\,x,\, b_{\,2},\, \cdots,\, b_{\,n}\,\right\| \,\leq\, 1\,\right\}$.
\item[$(II)$]\hspace{.2cm}$\|\,T\,\| \,=\, \sup\,\left\{\,\left|\,T\,(\,x,\, b_{\,2},\, \cdots,\, b_{\,n}\,)\,\right| \;:\; \left\|\,x,\, b_{\,2},\, \cdots,\, b_{\,n}\,\right\| \,=\, 1\,\right\}$.
\item[$(III)$]\hspace{.2cm}$ \|\,T\,\| \,=\, \sup\,\left \{\,\dfrac{\left|\,T\,(\,x,\, b_{\,2},\, \cdots,\, b_{\,n}\,)\,\right|}{\left\|\,x,\, b_{\,2},\, \cdots,\, b_{\,n}\,\right\|} \;:\; \left\|\,x,\, b_{\,2},\, \cdots,\, b_{\,n}\,\right\| \,\neq\, 0\,\right \}$. 
\end{itemize}
Also, we have \,$\left|\,T\,(\,x,\, b_{\,2},\, \cdots,\, b_{\,n}\,)\,\right| \,\leq\, \|\,T\,\|\, \left\|\,x,\, b_{\,2},\, \cdots,\, b_{\,n}\,\right\|\, \;\forall\; x \,\in\, W$.

In 1989, A.\,Misiak \cite{Misiak} introduced the concept of \,$n$-inner product for \,$n \,\geq\, 2$.

\begin{definition}\cite{Misiak}
Let \,$n \,\in\, \mathbb{N}$\; and \,$H$\, be a linear space of dimension greater than or equal to \,$n$\; over the field \,$\mathbb{K}$, where \,$\mathbb{K}$\, is the real or complex numbers field.\;An n-inner product on \,$H$\, is a map 
\[\left(\,x,\, y,\, x_{\,2},\, \cdots,\, x_{\,n}\,\right) \,\longmapsto\, \left<\,x,\, y \,|\, x_{\,2},\, \cdots,\, x_{\,n} \,\right>,\; x,\, y,\, x_{\,2},\, \cdots,\, x_{\,n} \,\in\, H\]from \,$H^{n \,+\, 1}$\, to the set \,$\mathbb{K}$\, such that for every \,$x,\, y,\, x_{\,1},\, x_{\,2},\, \cdots,\, x_{\,n} \,\in\, H$\, and \,$\alpha \,\in\, \mathbb{K}$,
\begin{itemize}
\item[(I)]\;\; $\left<\,x_{\,1},\, x_{\,1} \,|\, x_{\,2},\, \cdots,\, x_{\,n} \,\right> \,\geq\,  0$\; and \;$\left<\,x_{\,1},\, x_{\,1} \;|\; x_{\,2},\, \cdots,\, x_{\,n} \,\right> \;=\;  0$\; if and only if \;$x_{\,1},\, x_{\,2},\, \cdots,\, x_{\,n}$\; are linearly dependent,
\item[(II)]\;\; $\left<\,x,\, y \;|\; x_{\,2},\, \cdots,\, x_{\,n} \,\right> \;=\; \left<\,x,\, y \;|\; x_{\,i_{\,2}},\, \cdots,\, x_{\,i_{\,n}} \,\right> $\; for every permutations \\$\left(\, i_{\,2},\, \cdots,\, i_{\,n} \,\right)$\; of \;$\left(\, 2,\, \cdots,\, n \,\right)$,
\item[(III)]\;\; $\left<\,x,\, y \;|\; x_{\,2},\, \cdots,\, x_{\,n} \,\right> \;=\; \overline{\left<\,y,\, x \;|\; x_{\,2},\, \cdots,\, x_{\,n} \,\right> }$,
\item[(IV)]\;\; $\left<\,\alpha\,x,\, y \;|\; x_{\,2},\, \cdots,\, x_{\,n} \,\right> \;=\; \alpha \,\left<\,x,\, y \;|\; x_{\,2},\, \cdots,\, x_{\,n} \,\right> $,
\item[(V)]\;\; $\left<\,x \,+\, y,\, z \;|\; x_{\,2},\, \cdots,\, x_{\,n} \,\right> \;=\; \left<\,x,\, z \;|\; x_{\,2},\, \cdots,\, x_{\,n} \,\right> \,+\,  \left<\,y,\, z \;|\; x_{\,2},\, \cdots,\, x_{\,n} \,\right>$.
\end{itemize}
A linear space \,$H$\, together with an n-inner product \,$\left<\,\cdot,\, \cdot \,|\, \cdot,\, \cdots,\, \cdot\,\right>$\, is called an $n$-inner product space.
\end{definition}

\begin{theorem}\cite{Misiak}
For \,$n$-inner product space \,$\left(\,H,\, \left<\,\cdot,\, \cdot \,|\, \cdot,\, \cdots,\, \cdot\,\right>\,\right)$, 
\begin{equation}\label{0.en0.1}
\left|\,\left<\,x,\, y \,|\, x_{\,2},\,  \cdots,\, x_{\,n}\,\right>\,\right| \,\leq\, \left\|\,x,\, x_{\,2},\, \cdots,\, x_{\,n}\,\right\|\, \left\|\,y,\, x_{\,2},\, \cdots,\, x_{\,n}\,\right\|
\end{equation}
hold for all \,$x,\, y,\, x_{\,2},\, \cdots,\, x_{\,n} \,\in\, H$, where \[\left \|\,x_{\,1},\, x_{\,2},\, \cdots,\, x_{\,n}\,\right\| \,=\, \sqrt{\left <\,x_{\,1},\, x_{\,1} \;|\; x_{\,2},\,  \cdots,\, x_{\,n}\,\right>}\,.\] The inequality (\ref{0.en0.1}) is called Cauchy-Schwarz inequality.
\end{theorem}

An $n$-inner product space is called $n$-Hilbert space if it is complete with respect to it's induced \,$n$-norm.

\begin{definition}\cite{WR}
A complex algebra is a vector space \,$A$\, over the complex field \,$\mathbb{C}$\, in which a multiplication is defined that satisfies
\begin{align*}
x\,(\,y\,z\,) \,=\, (\,x\,y\,)\,z\,,\; \;(\,x \,+\, y\,)\,z \,=\, x\,z \,+\, y\,z\,,\; \;x\,(\,y \,+\, z\,) \,=\, x\,y \,+\, x\,z\,,
\end{align*}
and \,$\alpha\,(\,x\,y\,) \,=\, (\,\alpha\,x\,)\,y \,=\, x\,(\,\alpha\,y\,)$, for all \,$x,\, y$\, and \,$z$\, in \,$A$\, and for all scalars \,$\alpha$.
\end{definition}

\section{Best Approximation }

Let \,$X$\, be a \,$n$-normed linear space, \,$b_{\,2},\, b_{\,3},\, \cdots,\, b_{\,n}$\, be fixed elements in \,$X$\, and \,$G$\, be a subset of \,$X$.\,Suppose \,$x_{\,0} \,\in\, X$\, with \,$\left\{\,x_{\,0},\, b_{\,2},\, b_{\,3},\, \cdots,\, b_{\,n}\,\right\}$\, is linearly independent  and 
\[\delta_{\,b}\left(\,x_{\,0},\, G\,\right) \,=\, \inf\limits_{g \,\in\, G}\left\|\,x_{\,0} \,-\, g,\, b_{\,2},\, b_{\,3},\, \cdots,\, b_{\,n}\,\right\|.\]
Here \,$\delta_{\,b}\left(\,x_{\,0},\, G\,\right)$\, is the distance of \,$x_{\,0}$\, from \,$G$\, with respect to \,$b_{\,2},\, b_{\,3},\, \cdots,\, b_{\,n}$\, and sometimes it will be denoted by \,$\delta_{\,b}$.\,The problem for best approximation of \,$x_{\,0}$\, out of the elements of \,$G$\, is to find an element \,$g_{\,0} \,\in\, G$\, such that
\[\delta_{\,b} \,=\, \left\|\,x_{\,0} \,-\, g_{\,0},\, b_{\,2},\, b_{\,3},\, \cdots,\, b_{\,n}\,\right\|.\]
This element \,$g_{\,0} \,\in\, G$\, (\,if exists\,) has therefore minimum distance from \,$x_{\,0}$.\,Such an element \,$g_{\,0}$\, may or may not exist.\,We shall see below that under certain sufficient conditions, an element \,$g_{\,0}$\, with the above property may be found out.\,If \,$g_{\,0}$\, exists, then next problem arises about the uniqueness which has practical importance too.\,In this section, we shall deal with some of these problems.\,In the following theorem we observe that the finite dimensionality plays an important role for the existence of best approximation.

\begin{theorem}\label{3.thm3.101}
If \,$G$\, is a finite dimensional subspace of a \,$n$-normed linear space \,$X$, then for each \,$x$\, of \,$X$\, with \,$\left\{\,x,\, b_{\,2},\, b_{\,3},\, \cdots,\, b_{\,n}\,\right\}$\, is linearly independent, there exists a best approximation to \,$x$\, out of the elements of \,$G$.
\end{theorem}

\begin{proof}
For \,$x \,\in\, X$, consider the closed sphere 
\[\overline{B} \,=\, \left\{\,y \,\in\, G \,:\, \left\|\,y,\, b_{\,2},\, b_{\,3},\, \cdots,\, b_{\,n}\,\right\| \,\leq\, 2\,\left\|\,x ,\, b_{\,2},\, b_{\,3},\, \cdots,\, b_{\,n}\,\right\|\,\right\}.\]
Then \,$\theta \,\in\, \overline{B}$\, and for the distance of \,$x$\, from \,$\overline{B}$, we see that
\begin{align*}
\delta_{\,b}\,\left(\,x,\, \overline{B}\,\right) &\,=\, \inf\limits_{g \,\in\, \overline{B}}\,\left\|\,x \,-\, g,\, b_{\,2},\, b_{\,3},\, \cdots,\, b_{\,n}\,\right\| \\
&\,\leq\, \left\|\,x \,-\, \theta,\, b_{\,2},\, b_{\,3},\, \cdots,\, b_{\,n}\,\right\| \,=\, \left\|\,x,\, b_{\,2},\, b_{\,3},\, \cdots,\, b_{\,n}\,\right\|.
\end{align*}
If \,$y \,\notin\, \overline{B}$\, then 
\[\left\|\,y,\, b_{\,2},\, b_{\,3},\, \cdots,\, b_{\,n}\,\right\| \,>\, 2\,\left\|\,x ,\, b_{\,2},\, b_{\,3},\, \cdots,\, b_{\,n}\,\right\|\] 
and using the above 
\begin{align*}
\left\|\,x \,-\, y,\, b_{\,2},\, b_{\,3},\, \cdots,\, b_{\,n}\,\right\| &\,\geq\, \left\|\,y,\, b_{\,2},\, b_{\,3},\, \cdots,\, b_{\,n}\,\right\| \,-\, \left\|\,x,\, b_{\,2},\, b_{\,3},\, \cdots,\, b_{\,n}\,\right\|\\
&\,>\, \left\|\,x,\, b_{\,2},\, b_{\,3},\, \cdots,\, b_{\,n}\,\right\| \,\geq\, \delta_{\,b}\,\left(\,x,\, \overline{B}\,\right). 
\end{align*}
From this we see that \,$\delta_{\,b}\,\left(\,x,\, \overline{B}\,\right) \,=\, \delta_{\,b}\,\left(\,x,\, G\,\right)$\, and this value cannot be assumed by any \,$y \,\in\, G \,-\, \overline{B}$\, because of the strict inequality in the above relation.\,Therefore if a best approximation to \,$x$\, exists, it must be an element of \,$\overline{B}$.\,Since \,$\overline{B}$\, is closed and bounded subset of a finite dimensional subspace \,$G$, \,$\overline{B}$\, is compact.\,Consider the map 
\[f \,:\, y \,\mapsto\, \left\|\,x \,-\, y,\, b_{\,2},\, b_{\,3},\, \cdots,\, b_{\,n}\,\right\|\,,\, \,y \,\in\, \overline{B}.\]
Since the \,$n$-norm function \,$\left\|\,x \,-\, y,\, b_{\,2},\, b_{\,3},\, \cdots,\, b_{\,n}\,\right\|$\, (\,$y$\, variable\,) is continuous, \,$f$\, is a real valued continuous function and therefore \,$f$\, attains it's least upper bound.\,So, there is some \,$y_{\,0} \,\in\, \overline{B}$\, such that \,$f\,(\,y_{0}\,) \,=\, \inf\limits_{y \,\in\, \overline{B}}\,f\,(\,y\,)$.\,From this it follows that
\[\left\|\,x \,-\, y_{\,0},\, b_{\,2},\, b_{\,3},\, \cdots,\, b_{\,n}\,\right\| \,=\, \inf\limits_{y \,\in\, \overline{B}}\,\left\|\,x \,-\, y,\, b_{\,2},\, b_{\,3},\, \cdots,\, b_{\,n}\,\right\|.\] 
By definition, \,$y_{\,0}$\, is a best approximation to \,$x$\, out of the elements \,$G$.\,This proves the theorem.    
\end{proof}

In the next theorem, we see that the set of best approximations has certain geometric property.

\begin{theorem}
Let \,$G$\, be a subspace of a \,$n$-normed linear space \,$X$\, and \,$x \,\in\, X$\, with \,$\left\{\,x,\, b_{\,2},\, b_{\,3},\, \cdots,\, b_{\,n}\,\right\}$\, is linearly independent.\,Then the set \,$M$\, of best approximations to \,$x$\, out of \,$G$\, is convex.
\end{theorem} 

\begin{proof}
Let \,$\delta_{\,b} \,=\, \delta_{\,b}\,\left(\,x,\, G\,\right)$,\, i.\,e., \,$\delta_{\,b}$\, is the distance of \,$x_{\,0}$\, from \,$G$\, with respect to \,$b_{\,2},\, b_{\,3},\, \cdots,\, b_{\,n}$.\,If \,$M$\, is void or \,$M$\, contains only one element then there is nothing to prove.\,We therefore assume that \,$M$\, contains more than one element.\,Let \,$y,\, z \,\in\, M$,\, then from definition
\[\left\|\,x \,-\, y,\, b_{\,2},\, b_{\,3},\, \cdots,\, b_{\,n}\,\right\| \,=\, \delta_{\,b} \,,\; \; \left\|\,x \,-\, z,\, b_{\,2},\, b_{\,3},\, \cdots,\, b_{\,n}\,\right\| \,=\, \delta_{\,b}.\] 
Let \,$w$\, be any element of the form \,$w \,=\, \lambda\,y \,+\, \left(\,1 \,-\, \lambda\,\right)\,z$\, where \,$0 \,\leq\, \lambda \,\leq\, 1$.\,We should show that \,$w \,\in\, M$.\,Since \,$w \,\in\, G$,\, we see that \,$\left\|\,x \,-\, y,\, b_{\,2},\, b_{\,3},\, \cdots,\, b_{\,n}\,\right\| \,\geq\, \delta_{\,b}$\, and
\begin{align*}
&\left\|\,x \,-\, w,\, b_{\,2},\, b_{\,3},\, \cdots,\, b_{\,n}\,\right\|\\
& \,=\, \left\|\,\lambda\,\left(\,x \,-\, y\,\right) \,+\, \left(\,1 \,-\, \lambda\,\right)\,\left(\,x \,-\, z\,\right),\, b_{\,2},\, b_{\,3},\, \cdots,\, b_{\,n}\,\right\|\\
&\leq\,\lambda\,\left\|\,x \,-\, y,\, b_{\,2},\, b_{\,3},\, \cdots,\, b_{\,n}\,\right\| \,+\, \left(\,1 \,-\, \lambda\,\right)\,\left\|\,x \,-\, z,\, b_{\,2},\, b_{\,3},\, \cdots,\, b_{\,n}\,\right\|\\
&=\,\lambda\,\delta_{\,b} \,+\, \left(\,1 \,-\, \lambda\,\right)\,\delta_{\,b} \,=\, \delta_{\,b}\,,  
\end{align*}
so \,$\left\|\,x \,-\, w,\, b_{\,2},\, b_{\,3},\, \cdots,\, b_{\,n}\,\right\| \,=\, \delta_{\,b}$\, and \,$w \,\in\, M$.\,Since \,$y$\, and \,$z$\, are arbitrary elements of \,$M$,\, this proves the theorem.
\end{proof}     

Now, we exhibit some examples showing various possibilities of best approximations.

\begin{example}\label{3.examp.101}
\begin{itemize}
\item[$(i)$]Consider the space \,$C\,[\,a,\, b\,]$.\,For \,$x_{1} \,=\, x_{1}\,(\,t\,),\, x_{2} \,=\, x_{2}\,(\,t\,)$,\, \,$\cdots,\, x_{n} \,=\, x_{n}\,(\,t\,) \,\in\, C\,[\,a,\, b\,]$, define
\[\left\|\,x_{\,1},\, \cdots,\, x_{\,n}\,\right\|\]
\begin{equation}\label{eqpq2.1}
\,=\, \begin{cases}
\sup\limits_{a \,\leq\, t \,\leq\, b}\,\left|\,x_{1}\,(\,t\,)\,\right| \,\times\,\cdots\,\times\sup\limits_{a \,\leq\, t \,\leq\, b}\,\left|\,x_{n}\,(\,t\,)\,\right| & \text{if}\; x_{\,1},\, \cdots \,,\, x_{\,n} \;\text{are linearly }\\ &\;\hspace{.5cm}\text{independent,}\\0 & \text{if}\; x_{\,1},\, \cdots \,,\, x_{\,n} \;\text{are linearly}\\  &\;\hspace{.5cm}\text{dependent}. \end{cases}
\end{equation} 
Now, for every \,$x \,=\, x\,(\,t\,),\, y \,=\, y\,(\,t\,),\, x_{2} \,=\, x_{2}\,(\,t\,)$,\, \,$\cdots,\, x_{n} \,=\, x_{n}\,(\,t\,) \,\in\, C\,[\,a,\, b\,]$, and \,$\alpha \,\in\, \mathbb{R}$, we get
\begin{align*}
&\left\|\,x \,+\, y,\, x_{2},\,\cdots,\, x_{\,n}\,\right\|\\ 
&=\,\sup\limits_{a \,\leq\, t \,\leq\, b}\,\left|\,x\,(\,t\,) \,+\, y\,(\,t\,)\,\right| \,\times\, \sup\limits_{a \,\leq\, t \,\leq\, b}\,\left|\,x_{2}\,(\,t\,)\,\right|\,\times\, \cdots\,\times\sup\limits_{a \,\leq\, t \,\leq\, b}\,\left|\,x_{n}\,(\,t\,)\,\right|\\
&\leq\,\left(\,\sup\limits_{a \,\leq\, t \,\leq\, b}\,\left|\,x\,(\,t\,)\,\right| \,\times\, \sup\limits_{a \,\leq\, t \,\leq\, b}\,\left|\,x_{2}\,(\,t\,)\,\right|\,\times\, \cdots\,\times\sup\limits_{a \,\leq\, t \,\leq\, b}\,\left|\,x_{n}\,(\,t\,)\,\right|\,\right)\,+\\
&+\,\left(\,\sup\limits_{a \,\leq\, t \,\leq\, b}\,\left|\,y\,(\,t\,)\,\right| \,\times\, \sup\limits_{a \,\leq\, t \,\leq\, b}\,\left|\,x_{2}\,(\,t\,)\,\right|\,\times\, \cdots\,\times\sup\limits_{a \,\leq\, t \,\leq\, b}\,\left|\,x_{n}\,(\,t\,)\,\right|\,\right)\\
&=\,\left\|\,x,\, x_{2},\,\cdots,\, x_{\,n}\,\right\| \,+\, \left\|\,y,\, x_{2},\,\cdots,\, x_{\,n}\,\right\|
\end{align*}
and
\begin{align*}
&\left\|\,\alpha\,x,\, x_{\,2},\, \cdots,\, x_{\,n}\,\right\|\\
&=\,\sup\limits_{a \,\leq\, t \,\leq\, b}\,\left|\,\alpha\,x\,(\,t\,)\,\right| \,\times\, \sup\limits_{a \,\leq\, t \,\leq\, b}\,\left|\,x_{2}\,(\,t\,)\,\right|\,\times\, \cdots\,\times\sup\limits_{a \,\leq\, t \,\leq\, b}\,\left|\,x_{n}\,(\,t\,)\,\right|\\
&=\,|\,\alpha\,|\,\sup\limits_{a \,\leq\, t \,\leq\, b}\,\left|\,x\,(\,t\,)\,\right| \,\times\, \sup\limits_{a \,\leq\, t \,\leq\, b}\,\left|\,x_{2}\,(\,t\,)\,\right|\,\times\, \cdots\,\times\sup\limits_{a \,\leq\, t \,\leq\, b}\,\left|\,x_{n}\,(\,t\,)\,\right|\\
&=\,|\,\alpha\,|\,\left\|\,x,\, x_{\,2},\, \cdots,\, x_{\,n}\,\right\|.
\end{align*}
Therefore \,$C\,[\,a,\, b\,]$\, becomes a linear\;$n$-normed space with respect to the \,$n$-norm defined by (\ref{eqpq2.1}).\,Let \,\,$x_{k} \,=\, x_{k}\,(\,t\,) \,\in\, C\,[\,a,\, b\,]$\, be a Cauchy sequence, i.\,e, for every \,$e_{2} \,=\, e_{2}\,(\,t\,)$,\, \,$\cdots,\, e_{n} \,=\, e_{n}\,(\,t\,) \,\in\, C\,[\,a,\, b\,]$, \,$\left\|\,x_{k} \,-\, x_{l},\, e_{\,2},\, \cdots,\, e_{\,n}\,\right\| \,\to\, 0$\, as \,$k,\, l \,\to\, \infty$.\,Thus, for \,$\epsilon \,>\, 0$, there exists a positive integer \,$N$\, such that \,$\left\|\,x_{k} \,-\, x_{l},\, e_{\,2},\, \cdots,\, e_{\,n}\,\right\| \,<\, \epsilon$\, if \,$k,\, l \,\geq\, N$.\,This implies that if \,$k,\, l \,\geq\, N$,
\begin{align*}
&\sup\limits_{a \,\leq\, t \,\leq\, b}\,\left|\,x_{k}\,(\,t\,) \,-\, x_{l}\,(\,t\,)\,\right| \,\times\, \sup\limits_{a \,\leq\, t \,\leq\, b}\,\left|\,e_{2}\,(\,t\,)\,\right|\,\times\, \cdots\,\times\sup\limits_{a \,\leq\, t \,\leq\, b}\,\left|\,e_{n}\,(\,t\,)\,\right|  \,<\, \epsilon\\
&\Rightarrow\,\sup\limits_{a \,\leq\, t \,\leq\, b}\,\left|\,x_{k}\,(\,t\,) \,-\, x_{l}\,(\,t\,)\,\right| \,<\, \dfrac{\epsilon}{M^{\,n \,-\, 1}} \,=\, \epsilon^{\,\prime},
\end{align*}
where \,$M \,=\, \max\left\{\,\sup\limits_{a \,\leq\, t \,\leq\, b}\,\left|\,e_{2}\,(\,t\,)\,\right|,\, \cdots,\, \sup\limits_{a \,\leq\, t \,\leq\, b}\,\left|\,e_{n}\,(\,t\,)\,\right|\,\right\}$.\,Thus, we get that \,$\left|\,x_{k}\,(\,t\,) \,-\, x_{l}\,(\,t\,)\,\right| \,<\, \epsilon^{\,\prime}$\, for \,$k,\, l \,\geq\, N$\, and for all \,$t \,\in\, [\,a,\, b\,]$.\,This shows that the sequence \,$\left\{\,x_{k}\,(\,t\,)\,\right\}$\, converges uniformly to some function say, \,$x\,(\,t\,)$\, in \,$[\,a,\, b\,]$.\,Since \,$x_{k}\,(\,t\,)$\, are continuous, \,$x\,(\,t\,)$\, is also continuous on \,$[\,a,\, b\,]$\, and so \,\,$x \,\in\, C\,[\,a,\, b\,]$.\,Letting \,$l \,\to\, \infty$\, in the above inequality, we get  
\begin{align*}
&\sup\limits_{a \,\leq\, t \,\leq\, b}\,\left|\,x_{k}\,(\,t\,) \,-\, x\,(\,t\,)\,\right|\, \,\times\, \sup\limits_{a \,\leq\, t \,\leq\, b}\,\left|\,e_{2}\,(\,t\,)\,\right|\,\times\, \cdots\,\times\sup\limits_{a \,\leq\, t \,\leq\, b}\,\left|\,e_{n}\,(\,t\,)\,\right| \,<\, \epsilon,
\end{align*}
for \,$k \,\geq\, N$.\,This implies that 
\begin{align*}
\left\|\,x_{k} \,-\, x,\, e_{2},\,\cdots,\, e_{\,n}\,\right\|  \,<\, \epsilon,\; \;\text{for}\; \;k \,\geq\, N.
\end{align*}
So, \,$x_{k} \,\to\, x$\, as \,$k \,\to\, \infty$.\,This shows that \,$C\,[\,a,\, b\,]$\, is complete with respect to the above \,$n$-norm and hence it is a \,$n$-Banach space.

Now, let \,$x_{k}\,(\,t\,) \,=\, t^{\,k}$\, and \,$Y$\, be the subspace of \,$C\,[\,a,\, b\,]$\, generated by \,$x_{\,0},\, x_{\,1},\, \cdots,\, x_{\,m}$, where \,$m \,\geq\, n$\, is fixed.\,Then \,$Y$\, is finite dimensional and \,$Y$\, consists of all polynomials of degree at most \,$m$.\,If \,$x$\, is any element of \,$C\,[\,a,\, b\,]$\, then by Theorem \ref{3.thm3.101}, there exists a polynomial say, \,$p_{\,m}$\, of degree at most \,$m$\, such that
\begin{align*}
&\sup\limits_{a \,\leq\, t \,\leq\, b}\,\left|\,x\,(\,t\,) \,-\, p_{\,m}\,(\,t\,)\,\right|\, \,\times\, \sup\limits_{a \,\leq\, t \,\leq\, b}\,\left|\,e_{2}\,(\,t\,)\,\right|\,\times\, \cdots\,\times\sup\limits_{a \,\leq\, t \,\leq\, b}\,\left|\,e_{n}\,(\,t\,)\,\right|\\
&\leq\,\sup\limits_{a \,\leq\, t \,\leq\, b}\,\left|\,x\,(\,t\,) \,-\, p\,(\,t\,)\,\right|\, \,\times\, \sup\limits_{a \,\leq\, t \,\leq\, b}\,\left|\,e_{2}\,(\,t\,)\,\right|\,\times\, \cdots\,\times\sup\limits_{a \,\leq\, t \,\leq\, b}\,\left|\,e_{n}\,(\,t\,)\,\right|,
\end{align*}
for every \,$p\,(\,t\,) \,\in\, Y$. 
\item[$(ii)$]In this example, we shall see that the finite dimensionality of \,$G$\, in Theorem \ref{3.thm3.101} is essential.\,Consider the space \,$C\,[\,0,\, a\,]$\, where \,$0 \,<\, a \,<\, 1$\, and let \,$G$\, be the set of polynomials of any degree on \,$[\,0,\, a\,]$.\,Then the dimension of \,$G$\, infinite.\,Let \,$x\,(\,t\,) \,=\, (\,1 \,-\, t\,)^{\,-\, 1}$\, and \,$e_{j}\,(\,t\,) \,=\, (\,-\, 1\,)^{\,j}$, for \,$j \,=\, 2,\, \cdots,\, n$.\,Then for arbitrary \,$\epsilon \,>\, 0$\, there is \,$N$\, such that \,$\left\|\,x \,-\, y_{\,k},\, e_{2},\,\cdots,\, e_{\,n}\,\right\|  \,<\, \epsilon$\, for all \,$k \,\geq\, N$, where \,$y_{\,k}\,(\,t\,) \,=\, 1 \,+\, t \,+\, t^{\,2} \,+\, \cdots \,+\, t^{\,k} \,\in\, G$.\,Since \,$\epsilon \,>\, 0$\, is arbitrary, this shows that \,$\delta_{\,b}\,(\,x,\, G\,) \,=\, 0$.\,However, since \,$x$\, is not a polynomial, we see that there is no \,$y_{\,0} \,\in\, G$\, such that
\[\delta_{\,b}\,(\,x,\, G\,) \,=\, \left\|\,x \,-\, y_{\,0},\, e_{2},\,\cdots,\, e_{\,n}\,\right\| \,=\, 0.\]
\item[$(iii)$]In this example we show that the uniqueness of best approximations may not hold.\,Consider the space \,$\mathbb{R}^{\,n}$\, and for \,$x_{\,1} \,=\, \left(\,a^{(\,1\,)}_{\,1},\, a^{(\,1\,)}_{\,2},\, \cdots,\, a^{(\,1\,)}_{\,n}\,\right),\, x_{\,2} \,=\, \left(\,a^{(\,2\,)}_{\,1},\, a^{(\,2\,)}_{\,2},\, \cdots,\, a^{(\,2\,)}_{\,n}\,\right),\, \cdots,\,$\, \,$x_{\,n} \,=\, \left(\,a^{(\,n\,)}_{\,1},\, a^{(\,n\,)}_{\,2},\, \cdots,\, a^{(\,n\,)}_{\,n}\,\right) \,\in\, \mathbb{R}^{\,n}$\, define 
\[\left\|\,x_{\,1},\, x_{\,2},\, \cdots,\, x_{\,n}\,\right\|\]
\begin{equation}\label{eqpqe2.1}
\,=\, \begin{cases}
\max\limits_{\,i}\,\left|\,a^{(\,1\,)}_{\,i}\,\right| \,\times\, \max\limits_{\,i}\,\left|\,a^{(\,2\,)}_{\,i}\,\right| \,\times\, \,\cdots \,\times\, \max\limits_{\,i}\,\left|\,a^{(\,n\,)}_{\,i}\,\right| & \text{if}\; x_{\,1},\, \cdots \,,\, x_{\,n} \;\text{are linearly }\\ &\;\hspace{.5cm}\text{independent,}\\0 & \text{if}\; x_{\,1},\, \cdots \,,\, x_{\,n} \;\text{are linearly}\\  &\;\hspace{.5cm}\text{dependent}. \end{cases}
\end{equation}
Now, for every \,$x_{\,1} \,=\, \left(\,a^{(\,1\,)}_{\,1},\, a^{(\,1\,)}_{\,2},\, \cdots,\, a^{(\,1\,)}_{\,n}\,\right),\, y_{\,1} \,=\, \left(\,b^{(\,1\,)}_{\,1},\, b^{(\,1\,)}_{\,2},\, \cdots,\, b^{(\,1\,)}_{\,n}\,\right),\, x_{\,2} \,=\, \left(\,a^{(\,2\,)}_{\,1},\, a^{(\,2\,)}_{\,2},\, \cdots,\, a^{(\,2\,)}_{\,n}\,\right),\, \cdots,\,$\, \,$x_{\,n} \,=\, \left(\,a^{(\,n\,)}_{\,1},\, a^{(\,n\,)}_{\,2},\, \cdots,\, a^{(\,n\,)}_{\,n}\,\right) \,\in\, \mathbb{R}^{\,n}$\, and \,$\alpha \,\in\, \mathbb{R}$, we have 
\begin{align*}
&\left\|\,x_{\,1} \,+\, y_{\,1},\, x_{\,2},\, \cdots,\, x_{\,n}\,\right\|\\
&=\,\max\limits_{\,i}\,\left|\,a^{(\,1\,)}_{\,i} \,+\, b^{(\,1\,)}_{\,i}\,\right| \,\times\, \max\limits_{\,i}\,\left|\,a^{(\,2\,)}_{\,i}\,\right| \,\times\, \,\cdots \,\times\, \max\limits_{\,i}\,\left|\,a^{(\,n\,)}_{\,i}\,\right|\\
&\leq\,\left(\,\max\limits_{\,i}\,\left|\,a^{(\,1\,)}_{\,i}\,\right| \,\times\, \max\limits_{\,i}\,\left|\,a^{(\,2\,)}_{\,i}\,\right| \,\times\, \,\cdots \,\times\, \max\limits_{\,i}\,\left|\,a^{(\,n\,)}_{\,i}\,\right|\,\right)\,+\\
&\hspace{1cm}+\,\left(\,\max\limits_{\,i}\,\left|\,b^{(\,1\,)}_{\,i}\,\right| \,\times\, \max\limits_{\,i}\,\left|\,a^{(\,2\,)}_{\,i}\,\right| \,\times\, \,\cdots \,\times\, \max\limits_{\,i}\,\left|\,a^{(\,n\,)}_{\,i}\,\right|\,\right)\\
&=\,\left\|\,x_{\,1},\, x_{\,2},\, \cdots,\, x_{\,n}\,\right\| \,+\, \left\|\,y_{\,1},\, x_{\,2},\, \cdots,\, x_{\,n}\,\right\|,
\end{align*}
and 
\begin{align*}
&\left\|\,\alpha\,x_{\,1},\, x_{\,2},\, \cdots,\, x_{\,n}\,\right\|\\
&=\,\max\limits_{\,i}\,\left|\,\alpha\,a^{(\,1\,)}_{\,i}\,\right| \,\times\, \max\limits_{\,i}\,\left|\,a^{(\,2\,)}_{\,i}\,\right| \,\times\, \,\cdots \,\times\, \max\limits_{\,i}\,\left|\,a^{(\,n\,)}_{\,i}\,\right|\\
&=\,|\,\alpha\,|\,\left\|\,x_{\,1},\, x_{\,2},\, \cdots,\, x_{\,n}\,\right\|.
\end{align*}
Thus, \,$\mathbb{R}^{\,n}$\, is a linear \,$n$-normed space with respect to the above \,$n$-norm defined by (\ref{eqpqe2.1}).\,Let \,$G \,=\, \left\{\,y \,=\, (\,g,\, g,\, \cdots,\, g\,) \,:\, g \,\in\, \mathbb{R}\,\right\}$.\,Then \,$G$\, be a subspace of \,$\mathbb{R}^{\,n}$.\,Let \,$x \,=\, (\,1,\, 1,\, \cdots,\, 1\,)$\, and \,$e_{\,2} \,=\, \left(\,1,\, -\,1 ,\, \cdots,\, 1\,\right),\, \cdots,\,$\, \,$e_{\,n} \,=\, \left(\,1,\,  1,\, \cdots,\, \,-\, 1\,\right) \,\in\, \mathbb{R}^{\,n}$\,.\,Then for all \,$y \,\in\, G$, we have
\begin{align*}
&\left\|\,x \,-\, y,\, e_{\,2},\, \cdots,\, e_{\,n}\,\right\| \,=\, |\,1 \,-\, g\,| \,\geq\, 1.
\end{align*}
The distance from \,$x$\, to \,$G$\, with respect to \,$e_{\,2},\, \cdots,\, e_{\,n}$\, is \,$\delta_{\,b}\,\left(\,x,\, G\,\right) \,=\, 1$, and all \,$y \,=\, (\,g,\, g,\, \cdots,\, g\,)$\, with \,$|\,g\,| \,\leq\, 1$\, are best approximations to \,$x$\, out of \,$G$.   
\end{itemize}
\end{example}

\section{Strictly Convex $n$-Norm }

We observe from the above examples that there are spaces where the uniqueness of best approximations in \,$n$-normed space may not hold.\,To enquire about the uniqueness, the following definition becomes helpful in some cases.

\begin{definition}
A \,$n$-normed linear space \,$X$\, is called strictly convex if for all \,$x,\, y,\, e_{\,2},\, \cdots,\, e_{\,n} \,\in\, X$\, with \,$\left\|\,x,\, e_{\,2},\, \cdots,\, e_{\,n}\,\right\| \,=\, 1$,\, \,$\left\|\,y,\, e_{\,2},\, \cdots,\, e_{\,n}\,\right\| \,=\, 1$\, we have \,$\left\|\,x \,+\, y,\, e_{\,2},\, \cdots,\, e_{\,n}\,\right\| \,<\, 2$, unless \,$x \,=\, y$.\,A \,$n$-normed linear space with such a \,$n$-norm is called a strictly convex \,$n$-normed linear space.   
\end{definition}

If \,$\left\|\,x,\, e_{\,2},\, \cdots,\, e_{\,n}\,\right\| \,=\, 1$,\, \,$\left\|\,y,\, e_{\,2},\, \cdots,\, e_{\,n}\,\right\| \,=\, 1$\, then 
\[\left\|\,x \,+\, y,\, e_{\,2},\, \cdots,\, e_{\,n}\,\right\| \,\leq\, \left\|\,x,\, e_{\,2},\, \cdots,\, e_{\,n}\,\right\| \,+\, \left\|\,y,\, e_{\,2},\, \cdots,\, e_{\,n}\,\right\| \,=\, 2\]
 and the strict convexity excludes the equality sign, unless \,$x \,=\, y$.

In the following theorem, we observe that any \,$n$-Hilbert space is strictly convex.
\begin{theorem}
Any \,$n$-Hilbert space \,$H$\, is strictly convex.
\end{theorem}

\begin{proof}
Let \,$x,\, y \,\in\, H$, \,$x \,\neq\, y$\, and 
\[\left\|\,x,\, e_{\,2},\, \cdots,\, e_{\,n}\,\right\| \,=\, 1\,,\;\left\|\,y,\, e_{\,2},\, \cdots,\, e_{\,n}\,\right\| \,=\, 1.\] 
Let \,$\left\|\,x \,-\, y,\, e_{\,2},\, \cdots,\, e_{\,n}\,\right\| \,=\, \alpha$, then \,$\alpha \,>\, 0$.\,By the parallelogram law
\begin{align*}
&\left\|\,x \,+\, y,\, e_{\,2},\, \cdots,\, e_{\,n}\,\right\|^{\,2}\\
&=\,2\,\left(\,\left\|\,x,\, e_{\,2},\, \cdots,\, e_{\,n}\,\right\|^{\,2} \,+\, \left\|\,y,\, e_{\,2},\, \cdots,\, e_{\,n}\,\right\|^{\,2}\,\right) \,-\, \left\|\,x \,-\, y,\, e_{\,2},\, \cdots,\, e_{\,n}\,\right\|^{\,2}\\
&=\,2\,(\,1 \,+\, 1\,) \,-\, \alpha^{\,2} \,=\, 4 \,-\, \alpha^{\,2} \,<\, 4
\end{align*}
and this gives \,$\left\|\,x \,+\, y,\, e_{\,2},\, \cdots,\, e_{\,n}\,\right\| \,<\, 2$.
\end{proof}

\begin{example}
We consider the space \,$C\,[\,a,\, b\,]$\, with respect to the \,$n$-norm given by (\ref{eqpq2.1}).\,Let \,$x,\, y,\, e_{2},\, \cdots,\, e_{n}$\, be defined by
\[x\,(\,t\,) \,=\, 1\,,\, y\,(\,t\,) \,=\, \dfrac{t \,-\, a}{b \,-\, a}\,,\, e_{k}\,(\,t\,) \,=\, (\,-\, 1\,)^{\,k}\,,\, \,k \,=\, 2,\, \cdots,\, n\,,\]
where \,$t \,\in\, [\,a,\, b\,]$.\,Then \,$x \,\neq\, y$\, and \,$x,\, y,\, e_{2},\, \cdots,\, e_{n} \,\in\, C\,[\,a,\, b\,]$.\,Also,
\[\left\|\,x,\, e_{\,2},\, \cdots,\, e_{\,n}\,\right\| \,=\, 1\,, \,\left\|\,y,\, e_{\,2},\, \cdots,\, e_{\,n}\,\right\| \,=\, 1.\]
Now,
\begin{align*}
&\left\|\,x \,+\, y,\, e_{\,2},\, \cdots,\, e_{\,n}\,\right\| \\
&=\,\sup\limits_{a \,\leq\, t \,\leq\, b}\,\left|\,1 \,+\, \dfrac{t \,-\, a}{b \,-\, a}\,\right| \,\times\,\sup\limits_{a \,\leq\, t \,\leq\, b}\,\left|\,e_{2}\,(\,t\,)\,\right|\,\times\cdots\,\times\sup\limits_{a \,\leq\, t \,\leq\, b}\,\left|\,e_{n}\,(\,t\,)\,\right|\\
&=\,2.
\end{align*} 
Therefore \,$C\,[\,a,\, b\,]$\, is not strictly convex. 
\end{example}

In the next theorem, we see that strict convexity implies the uniqueness of best approximations. 

\begin{theorem}
Let \,$X$\, be a \,$n$-normed linear space with a strictly convex \,$n$-norm and \,$G$\, be a subspace of \,$X$.\,If \,$x \,\in\, X$\, then there is at most one best approximation to \,$x$\, out of the elements of \,$G$.
\end{theorem}

\begin{proof}
Let  \,$\delta_{\,b} \,=\, \textit{dist}\,\left(\,x,\, G\,\right)$\, and suppose, if possible that there are two distinct best approximations \,$g_{\,1}$\, and \,$g_{\,2}$\, to \,$x$.\,Then
\[\left\|\,x \,-\, g_{\,1},\, b_{\,2},\, \cdots,\, b_{\,n}\,\right\| \,=\, \delta_{\,b}\,,\;\left\|\,x \,-\, g_{\,2},\, b_{\,2},\, \cdots,\, b_{\,n}\,\right\| \,=\, \delta_{\,b}\]
and \,$x \,-\, g_{\,1}$\, and \,$x \,-\, g_{\,2}$\, are distinct.\,Therefore by strict convexity
\[\left\|\,x \,-\, g_{\,1} \,+\, x \,-\, g_{\,2},\, b_{\,2},\, \cdots,\, b_{\,n}\,\right\| \,<\, 2\,\delta_{\,b}\]and this implies that 
\[\left\|\,x \,-\, \dfrac{\left(\,g_{\,1} \,+\, g_{\,2}\,\right)}{2},\, b_{\,2},\, \cdots,\, b_{\,n}\,\right\| \,<\, \delta_{\,b}.\]Because \,$\dfrac{\left(\,g_{\,1} \,+\, g_{\,2}\,\right)}{2} \,\in\, G$, this implies a contradiction and the theorem is proved. 
\end{proof}

\begin{theorem}\label{th1.11}
Let \,$C$\, be a closed convex subset of a \,$n$-Hilbert space \,$H$.\,Then \,$C$\, contains a unique element of smallest \,$n$-norm.
\end{theorem}

\begin{proof}
Let \,$d \,=\, \inf\,\left\{\,\left\|\,x,\, b_{\,2},\, \cdots,\, b_{\,n}\,\right\| \,:\, x \,\in\, C\,\right\}$.\,Then there exists a sequence \,$\left\{\,x_{\,k}\,\right\}$\, of elements of \,$C$\, such that \,$\left\|\,x_{\,k},\, b_{\,2},\, \cdots,\, b_{\,n}\,\right\| \,\to\, d$.\,Now, \,$x_{\,k},\, x_{\,l} \,\in\, C$\, and \,$C$\, is convex, so \,$\dfrac{x_{\,k} \,+\, x_{\,l}}{2} \,\in\, C$\, and therefore
\[\left\|\,\dfrac{x_{\,k} \,+\, x_{\,l}}{2},\, b_{\,2},\, \cdots,\, b_{\,n}\,\right\| \,\geq\, d.\]This implies that
\[\left\|\,x_{\,k} \,+\, x_{\,l},\, b_{\,2},\, \cdots,\, b_{\,n}\,\right\| \,\geq\, 2\,d.\]Using parallelogram law, we get
\begin{align}
&\left\|\,x_{\,k} \,-\, x_{\,l},\, b_{\,2},\, \cdots,\, b_{\,n}\,\right\|^{\,2}\nonumber\\
&=\,2\,\left\|\,x_{\,k},\, b_{\,2},\, \cdots,\, b_{\,n}\,\right\|^{\,2} \,+\, 2\,\left\|\,x_{\,l},\, b_{\,2},\, \cdots,\, b_{\,n}\,\right\|^{\,2} \,-\, \left\|\,x_{\,k} \,+\, x_{\,l},\, b_{\,2},\, \cdots,\, b_{\,n}\,\right\|^{\,2}\nonumber \\
&\leq\,2\,\left\|\,x_{\,k},\, b_{\,2},\, \cdots,\, b_{\,n}\,\right\|^{\,2} \,+\, 2\,\left\|\,x_{\,l},\, b_{\,2},\, \cdots,\, b_{\,n}\,\right\|^{\,2} \,-\, 4\,d^{\,2}.\label{eq1.11}
\end{align}
Since \,$\left\|\,x_{\,k},\, b_{\,2},\, \cdots,\, b_{\,n}\,\right\| \,\to\, d$\, as \,$k \,\to\, \infty$\, and \,$\left\|\,x_{\,l},\, b_{\,2},\, \cdots,\, b_{\,n}\,\right\| \,\to\, d$\, as \,$l \,\to\, \infty$,\, the right hand side of (\ref{eq1.11}) tends to \,$2\,d^{\,2} \,+\, 2\,d^{\,2} \,-\, 4\,d^{\,2} \,=\, 0$\, as \,$k,\, l \,\to\, \infty$. So, \,$\left\{\,x_{\,k}\,\right\}$\, is a Cauchy sequence.\,Since \,$H$\, is complete, there exists \,$x$\, such that \,$x_{\,k} \,\to\, x$.\,As \,$C$\, is closed, we have \,$x \,\in\, C$.\,So,
\[\left\|\,x,\, b_{\,2},\, \cdots,\, b_{\,n}\,\right\| \,=\, \lim\limits_{k \,\to\, \infty}\,\left\|\,x_{\,k},\, b_{\,2},\, \cdots,\, b_{\,n}\,\right\| \,=\, d.\]Therefore, \,$x$\, is an element in \,$C$\, with the smallest \,$n$-norm (\,i.\,e.,\, $d$\,).\,Now, we prove the uniqueness.\,Let \,$y$\, be an element of \,$C$\, such that \,$\left\|\,y,\, b_{\,2},\, \cdots,\, b_{\,n}\,\right\|$\, also equals \,$d$.\,Then \,$C$\, being convex, we have
\,$\dfrac{x \,+\, y}{2} \,\in\, C$\, and also
\begin{equation}\label{eq1.12}
\left\|\,\dfrac{x \,+\, y}{2},\, b_{\,2},\, \cdots,\, b_{\,n}\,\right\| \,\geq\, d.
\end{equation}
Applying the parallelogram law, we obtain
\begin{align*}
&\left\|\,\dfrac{x \,+\, y}{2},\, b_{\,2},\, \cdots,\, b_{\,n}\,\right\|^{\,2}\\
&=\,\dfrac{\left\|\,x,\, b_{\,2},\, \cdots,\, b_{\,n}\,\right\|^{\,2}}{2} \,+\, \dfrac{\left\|\,y,\, b_{\,2},\, \cdots,\, b_{\,n}\,\right\|^{\,2}}{2} \,-\, \dfrac{\left\|\,x \,-\, y,\, b_{\,2},\, \cdots,\, b_{\,n}\,\right\|^{\,2}}{2}\\
&<\,\dfrac{\left\|\,x,\, b_{\,2},\, \cdots,\, b_{\,n}\,\right\|^{\,2}}{2} \,+\, \dfrac{\left\|\,y,\, b_{\,2},\, \cdots,\, b_{\,n}\,\right\|^{\,2}}{2} \,=\, d^{\,2}\; \;\;\textit{if}\;\; \;x \,\neq\, y
\end{align*}
and this gives that 
\[\left\|\,\dfrac{x \,+\, y}{2},\, b_{\,2},\, \cdots,\, b_{\,n}\,\right\| \,<\, d.\]\,This contradicts (\ref{eq1.12}) and so \,$x \,=\, y$.\,This proves the theorem.   
\end{proof}

\begin{theorem}\label{th1.12}
Let \,$M$\, be a closed subspace of a \,$n$-Hilbert space \,$H$\, and \,$x \,\in\, H \,-\, M$.\,Let \,$d$\, be the distance from \,$x$\, to \,$M$\, with respect to \,$b_{\,2},\, \cdots,\, b_{\,n}$,\, i.\,e.,
\[d \,=\, \inf\limits_{y \,\in\, M}\left\|\,x \,-\, y,\, b_{\,2},\, b_{\,3},\, \cdots,\, b_{\,n}\,\right\|.\]
Then there exists a unique element \,$y_{\,0} \,\in\, M$\, such that 
\[\left\|\,x \,-\, y_{\,0},\, b_{\,2},\, b_{\,3},\, \cdots,\, b_{\,n}\,\right\| \,=\, d.\]
\end{theorem}

\begin{proof}
Consider the set \,$C$\, defined by
\[C \,=\, \left\{\,x \,+\, y \,:\, y \,\in\, M\,\right\}.\]
We will now see that \,$C$\, is a closed convex set and that \,$d$\, is the distance from the origin \,$\theta$\, to \,$C$.
\begin{itemize}
\item[$(1)$]\,$C$\, is closed:\hspace{.4cm}Let \,$x \,+\, y_{\,k} \,\to\, x \,+\, y$\, where \,$y_{\,k} \,\in\, M$.\,So,
\[\left\|\,x \,+\, y_{\,k} \,-\, x \,-\, y,\, b_{\,2},\, b_{\,3},\, \cdots,\, b_{\,n}\,\right\| \,=\, \left\|\,y_{\,k} \,-\, y,\, b_{\,2},\, b_{\,3},\, \cdots,\, b_{\,n}\,\right\| \,\to\, 0\]
i.\,e., \,$y_{\,k} \,\to\, y$.\,Since \,$M$\, is closed, \,$y \,\in\, M$\, and this implies that \,$x \,+\, y \,\in\, C$.
\item[$(2)$]\,$C$\, is convex:\hspace{.4cm}Let \,$z_{\,1},\, z_{\,2} \,\in\, C$, then \,$z_{\,1} \,=\, x \,+\, y_{\,1}$,\, \,$z_{\,2} \,=\, x \,+\, y_{\,2}$, where \,$y_{\,1},\, y_{\,2} \,\in\, M$.\,So, if \,$0 \,\leq\, \lambda \,\leq\, 1$, then 
\begin{align*}
\lambda\,z_{\,1} \,+\, \left(\,1 \,-\, \lambda\,\right)\,z_{\,2} \,=\, x \,+\,\lambda\,y_{\,1} \,+\, \left(\,1 \,-\, \lambda\,\right)\,y_{\,2} \,=\, x \,+\, y\,, 
\end{align*}
where \,$y \,=\, \lambda\,y_{\,1} \,+\, \left(\,1 \,-\, \lambda\,\right)\,y_{\,2} \,\in\, M$, because \,$M$\, is a subspace.\,So, \,$\lambda\,z_{\,1} \,+\, \left(\,1 \,-\, \lambda\,\right)\,z_{\,2} \,\in\, C$\, and \,$C$\, is convex.
\item[$(3)$]\,$d$\, is the distance from the origin to \,$C$:\hspace{.2cm} We have
\begin{align*}
d &\,=\, \inf\limits_{y^{\,\prime} \,\in\, M}\left\|\,x \,-\, y^{\,\prime},\, b_{\,2},\, b_{\,3},\, \cdots,\, b_{\,n}\,\right\|\\
&=\,\inf\limits_{y \;=\, \,-\, y^{\,\prime} \,\in\, M}\left\|\,x \,+\, y,\, b_{\,2},\, b_{\,3},\, \cdots,\, b_{\,n}\,\right\|\\ 
&=\,\inf\limits_{z \,=\, x \,+\, y  \,\in\, C}\left\|\,z,\, b_{\,2},\, b_{\,3},\, \cdots,\, b_{\,n}\,\right\|.
\end{align*}
By Theorem \ref{th1.11}, there exists a unique element \,$z \,\in\, C$\, such that 
\[\left\|\,z,\, b_{\,2},\, b_{\,3},\, \cdots,\, b_{\,n}\,\right\| \,=\, d.\]
Let \,$y_{\,0} \,=\, x \,-\, z$\, then clearly \,$y_{\,0} \,\in\, M$\, and 
\[\left\|\,x \,-\, y_{\,0},\, b_{\,2},\, b_{\,3},\, \cdots,\, b_{\,n}\,\right\| \,=\, \left\|\,z,\, b_{\,2},\, b_{\,3},\, \cdots,\, b_{\,n}\,\right\| \,=\, d.\]
We now prove the uniqueness of \,$y_{\,0}$.

Let \,$y_{\,1} \,\in\, M$\, and \,$y_{\,1} \,\neq\, y_{\,0}$\, and if possible, suppose that 
\[\left\|\,x \,-\, y_{\,1},\, b_{\,2},\, b_{\,3},\, \cdots,\, b_{\,n}\,\right\| \,=\, d.\]Then \,$z_{\,1} \,=\, x \,-\, y_{\,1} \,\in\, C$\, such that \,$z_{\,1} \,\neq\, z_{\,0}$\, and moreover \[\left\|\,z_{\,1},\, b_{\,2},\, b_{\,3},\, \cdots,\, b_{\,n}\,\right\| \,=\, d.\]This contradicts the uniqueness of \,$z_{\,0}$.\,Hence the theorem is proved.
\end{itemize}
\end{proof}

In Theorem \ref{th1.12}, if instead of assuming \,$M$\, to be a closed subspace, we assume \,$M$\, to be a closed convex subset of \,$H$\,, we still obtain the best approximation as the next theorem shows.

\begin{theorem}
Let \,$M$\, be a nonempty closed convex subset of \,$H$\, and \,$x \,\in\, H $.\,Let \,$d$\, be the distance from \,$x$\, to \,$M$\, with respect to \,$b_{\,2},\, \cdots,\, b_{\,n}$,\, i.\,e.,
\[d \,=\, \inf\limits_{y \,\in\, M}\left\|\,x \,-\, y,\, b_{\,2},\, b_{\,3},\, \cdots,\, b_{\,n}\,\right\|.\]
Then there exists a unique element \,$y_{\,0} \,\in\, M$\, such that 
\[\left\|\,x \,-\, y_{\,0},\, b_{\,2},\, b_{\,3},\, \cdots,\, b_{\,n}\,\right\| \,=\, d.\]
\end{theorem}

\begin{proof}
It follows from the definition \,$d$\, that there exists a sequence \,$\left\{\,y_{\,k}\,\right\}$\, in \,$M$\, such that 
\begin{equation}\label{eq1.13}
d \,=\, \lim\limits_{k \,\to\, \infty}\left\|\,x \,-\, y_{\,k},\, b_{\,2},\, b_{\,3},\, \cdots,\, b_{\,n}\,\right\|.
\end{equation}
Using parallelogram law, we see that 
\begin{align*}
&\left\|\,y_{\,k} \,-\, y_{\,l},\, b_{\,2},\, b_{\,3},\, \cdots,\, b_{\,n}\,\right\|^{\,2}\\
&=\,2\,\left\|\,y_{\,k} \,-\, x,\, b_{\,2},\, b_{\,3},\, \cdots,\, b_{\,n}\,\right\|^{\,2} \,+\, 2\,\left\|\,y_{\,l} \,-\, x,\, b_{\,2},\, b_{\,3},\, \cdots,\, b_{\,n}\,\right\|^{\,2}\,-\\
&\hspace{1.3cm}\,-\left\|\,2\,x \,-\, y_{\,k} \,-\, y_{\,l},\, b_{\,2},\, b_{\,3},\, \cdots,\, b_{\,n}\,\right\|^{\,2}\\
&=\,2\,\left\|\,y_{\,k} \,-\, x,\, b_{\,2},\, b_{\,3},\, \cdots,\, b_{\,n}\,\right\|^{\,2} \,+\, 2\,\left\|\,y_{\,l} \,-\, x,\, b_{\,2},\, b_{\,3},\, \cdots,\, b_{\,n}\,\right\|^{\,2}\,-\\
&\hspace{1.3cm}\,-\, 4\,\left\|\,x \,-\, \dfrac{1}{2}\,\left(\,y_{\,k} \,+\, y_{\,l}\,\right),\, b_{\,2},\, b_{\,3},\, \cdots,\, b_{\,n}\,\right\|^{\,2}.
\end{align*}
As \,$M$\, is convex, \,$\dfrac{1}{2}\,\left(\,y_{\,k} \,+\, y_{\,l}\,\right) \,\in\, M$\, and therefore 
\[\left\|\,x \,-\, \dfrac{1}{2}\,\left(\,y_{\,k} \,+\, y_{\,l}\,\right),\, b_{\,2},\, b_{\,3},\, \cdots,\, b_{\,n}\,\right\|^{\,2} \,\geq\, d^{\,2}.\]
Using this we obtain from the above
\begin{align*}
&\left\|\,y_{\,k} \,-\, y_{\,l},\, b_{\,2},\, b_{\,3},\, \cdots,\, b_{\,n}\,\right\|^{\,2}\\
& \,\leq\, 2\,\left\|\,y_{\,k} \,-\, x,\, b_{\,2},\, b_{\,3},\, \cdots,\, b_{\,n}\,\right\|^{\,2} \,+\, 2\,\left\|\,y_{\,l} \,-\, x,\, b_{\,2},\, b_{\,3},\, \cdots,\, b_{\,n}\,\right\|^{\,2} \,-\, 4\,d^{\,2}
\end{align*} 
and now from (\ref{eq1.13}) it follows that
\[\left\|\,y_{\,k} \,-\, y_{\,l},\, b_{\,2},\, b_{\,3},\, \cdots,\, b_{\,n}\,\right\| \,\to\, \infty\; \;\;\textit{as}\;\; \;k,\, l \,\to\, \infty\]
i.\,e., \,$\left\{\,y_{\,k}\,\right\}$\, is a Cauchy sequence.\,Consequently there exists \,$y_{\,0} \,\in\, H$\, such that \,$\left\|\,y_{\,k} \,-\, y_{\,0},\, b_{\,2},\, b_{\,3},\, \cdots,\, b_{\,n}\,\right\| \,\to\, \infty$\, as \,$k \,\to\, \infty$.\,Since \,$M$\, is closed, \,$y_{\,0}$\, also belongs to \,$M$.\,Now,
\begin{align*}
&\left\|\,x \,-\, y_{\,0},\, b_{\,2},\, b_{\,3},\, \cdots,\, b_{\,n}\,\right\|\\
&\leq\,\left\|\,x \,-\, y_{\,k},\, b_{\,2},\, b_{\,3},\, \cdots,\, b_{\,n}\,\right\| \,+\, \left\|\,y_{\,k} \,-\, y_{\,0},\, b_{\,2},\, b_{\,3},\, \cdots,\, b_{\,n}\,\right\|
\end{align*}
and because 
\[\left\|\,x \,-\, y_{\,k},\, b_{\,2},\, b_{\,3},\, \cdots,\, b_{\,n}\,\right\| \,\to\, d\,,\; \;\left\|\,y_{\,k} \,-\, y_{\,0},\, b_{\,2},\, b_{\,3},\, \cdots,\, b_{\,n}\,\right\| \,\to\, 0\,,\]
it follows that \,$\left\|\,x \,-\, y_{\,0},\, b_{\,2},\, b_{\,3},\, \cdots,\, b_{\,n}\,\right\| \,\leq\, d$.\,Also, it is easy to clear that \\$\left\|\,x \,-\, y_{\,0},\, b_{\,2},\, b_{\,3},\, \cdots,\, b_{\,n}\,\right\| \,\geq\, d$.\,Thus, \,$\left\|\,x \,-\, y_{\,0},\, b_{\,2},\, b_{\,3},\, \cdots,\, b_{\,n}\,\right\| \,=\, d$.

We now prove the uniqueness part.\,Assume \,$y_{\,1} \,\in\, M$\, be such that \\$\left\|\,x \,-\, y_{\,1},\, b_{\,2},\, b_{\,3},\, \cdots,\, b_{\,n}\,\right\| \,=\, d$.\,Convexity of \,$M$\, implies that \,$\dfrac{1}{2}\,\left(\,y_{\,0} \,+\, y_{\,1}\,\right) \,\in\, M$.\,From these we obtain that
\begin{align*}
d &\,\leq\, \left\|\,x \,-\, \dfrac{1}{2}\,\left(\,y_{\,0} \,+\, y_{\,1}\,\right),\, b_{\,2},\, b_{\,3},\, \cdots,\, b_{\,n}\,\right\|\\
& =\,\left\|\,\dfrac{1}{2}\,x \,-\, \dfrac{1}{2}\,y_{\,0} \,+\, \dfrac{1}{2}\,x \,-\, \dfrac{1}{2}\,y_{\,1},\, b_{\,2},\, b_{\,3},\, \cdots,\, b_{\,n}\,\right\|\\
&\leq\,\dfrac{1}{2}\,\left\|\,x \,-\, y_{\,0},\, b_{\,2},\, b_{\,3},\, \cdots,\, b_{\,n}\,\right\| \,+\, \dfrac{1}{2}\,\left\|\,x \,-\, y_{\,1},\, b_{\,2},\, b_{\,3},\, \cdots,\, b_{\,n}\,\right\| \,=\, d
\end{align*} 
and therefore 
\[\left\|\,x \,-\, \dfrac{1}{2}\,\left(\,y_{\,0} \,+\, y_{\,1}\,\right),\, b_{\,2},\, b_{\,3},\, \cdots,\, b_{\,n}\,\right\| \,=\, d.\]
Applying parallelogram law we now see that
\begin{align*}
&\left\|\,y_{\,0} \,-\, y_{\,1},\, b_{\,2},\, b_{\,3},\, \cdots,\, b_{\,n}\,\right\|^{\,2}\\
&=\,2\,\left\|\,x \,-\, y_{\,0},\, b_{\,2},\, b_{\,3},\, \cdots,\, b_{\,n}\,\right\|^{\,2} \,+\, 2\,\left\|\,x \,-\, y_{\,1},\, b_{\,2},\, b_{\,3},\, \cdots,\, b_{\,n}\,\right\|^{\,2}\,-\\
&\hspace{1.3cm}\,-\, 4\,\left\|\,x \,-\, \dfrac{1}{2}\,\left(\,y_{\,0} \,+\, y_{\,1}\,\right),\, b_{\,2},\, b_{\,3},\, \cdots,\, b_{\,n}\,\right\|^{\,2}\\
&=\,2\,d^{\,2} \,+\, 2\,d^{\,2} \,-\, 4\,d^{\,2} \,=\, 0,
\end{align*}
i.\,e., \,$y_{\,0} \,=\, y_{\,1}$.\,This proves the theorem.
\end{proof}

\begin{definition}
The \,$n$-Banach space \,$X$\, is said to be uniformly convex if for every \,$\epsilon \,>\, 0$\, there exists \,$\delta \,>\, 0$\, such that if 
\[\left\|\,x,\, b_{\,2},\, b_{\,3},\, \cdots,\, b_{\,n}\,\right\| \,\leq\, 1\,,\;\; \left\|\,y,\, b_{\,2},\, b_{\,3},\, \cdots,\, b_{\,n}\,\right\| \,\leq\, 1\]
and 
\[\left\|\,x \,-\, y,\, b_{\,2},\, b_{\,3},\, \cdots,\, b_{\,n}\,\right\| \,\geq\, \epsilon\]
then
\[\left\|\,\dfrac{1}{2}\,\left(\,x \,+\, y\,\right),\, b_{\,2},\, b_{\,3},\, \cdots,\, b_{\,n}\,\right\| \,\leq\, 1 \,-\, \delta\]
\end{definition}
 
\begin{theorem}
Every \,$n$-inner product space \,$H$\, is uniformly convex with respect to the \,$n$-norm defined by the \,$n$-inner product.\,In particular, every \,$n$-Hilbert space is uniformly convex.
\end{theorem} 
 
\begin{proof}
Let \,$\epsilon \,>\, 0$\, be arbitrary and 
\[\left\|\,x,\, b_{\,2},\, b_{\,3},\, \cdots,\, b_{\,n}\,\right\| \,\leq\, 1\,,\;\; \left\|\,y,\, b_{\,2},\, b_{\,3},\, \cdots,\, b_{\,n}\,\right\| \,\leq\, 1\,,\]
\[\left\|\,x \,-\, y,\, b_{\,2},\, b_{\,3},\, \cdots,\, b_{\,n}\,\right\| \,\geq\, \epsilon.\]
By the parallelogram law we obtain
\begin{align*}
&\left\|\,\dfrac{1}{2}\,\left(\,x \,+\, y\,\right),\, b_{\,2},\, b_{\,3},\, \cdots,\, b_{\,n}\,\right\|^{\,2}\\
&=\,2\,\left(\,\left\|\,\dfrac{x}{2},\, b_{\,2},\, \cdots,\, b_{\,n}\,\right\|^{\,2} \,+\, \left\|\,\dfrac{y}{2},\, b_{\,2},\, \cdots,\, b_{\,n}\,\right\|^{\,2}\,\right) \,-\, \left\|\,\dfrac{x \,-\, y}{2},\, b_{\,2},\, \cdots,\, b_{\,n}\,\right\|^{\,2}\\
&\leq\,2\,\left(\,\dfrac{1}{4} \,+\, \dfrac{1}{4}\,\right) \,-\, \left(\,\dfrac{\epsilon}{2}\,\right)^{\,2} \,=\, 1 \,-\, \dfrac{\epsilon^{\,2}}{4}\\
&\Rightarrow\, \left\|\,\dfrac{1}{2}\,\left(\,x \,+\, y\,\right),\, b_{\,2},\, b_{\,3},\, \cdots,\, b_{\,n}\,\right\| \,\leq\, \left(\,1 \,-\, \dfrac{\epsilon^{\,2}}{4}\,\right)^{1 \,/\, 2}
\end{align*}
and this gives that \,$\delta\,(\,\epsilon\,) \,\geq\, 1 \,-\,  \left(\,1 \,-\, \dfrac{\epsilon^{\,2}}{4}\,\right)^{1 \,/\, 2}$\, which is positive for all positive \,$\epsilon\,(\,\leq\, 2\,)$.\,Hence \,$H$\, is uniformly convex.\,This proves the theorem.
\end{proof}

\begin{theorem}
Every uniformly convex \,$n$-Banach space is strictly convex.
\end{theorem}

\begin{proof}
This follows from the definition.
\end{proof}

\begin{theorem}
Let \,$X$\, be a uniformly convex \,$n$-Banach space and \,$\left\{\,x_{\,k}\,\right\}$\, be a sequence in \,$X$\, such that \,$\left\|\,x_{\,k},\, b_{\,2},\, b_{\,3},\, \cdots,\, b_{\,n}\,\right\| \,\to\, 1$\, as \,$k \,\to\, \infty$\, and 
\[\left\|\,x_{\,k} \,+\, x_{\,l},\, b_{\,2},\, b_{\,3},\, \cdots,\, b_{\,n}\,\right\| \,\to\, 2\; \;\;\textit{as}\; \;k,\, l \,\to\, \infty.\]\,Then \,$\left\{\,x_{\,k}\,\right\}$\, is a Cauchy sequence. 
\end{theorem} 
 
\begin{proof}
If \,$\left\{\,x_{\,k}\,\right\}$\, is not a Cauchy sequence then there exists \,$\epsilon \,>\, 0$\, and a sequence of positive integers \,$k_{\,1} \,<\, k_{\,2} \,<\, \cdots$\, such that 
\[\left\|\,x_{\,k_{\,j}},\, b_{\,2},\, b_{\,3},\, \cdots,\, b_{\,n}\,\right\| \,<\, 1 \,+\, \dfrac{1}{j}\,,\; \left\|\,x_{\,k_{\,j \,+\, 1}} \,-\, x_{\,k_{\,j}},\, b_{\,2},\, b_{\,3},\, \cdots,\, b_{\,n}\,\right\| \,\geq\, 2\,\epsilon\]
for \,$j \,=\, 1,\, 2,\, \cdots$.\,Since \,$X$\, is uniformly convex, corresponding to this \,$\epsilon$,\, there exists \,$\delta \,>\, 0$\, such that
\[\left\|\,\dfrac{1}{2}\,\left(\,x \,+\, y\,\right),\, b_{\,2},\, b_{\,3},\, \cdots,\, b_{\,n}\,\right\| \,\leq\, 1 \,-\, \delta\]
whenever
\[\left\|\,x,\, b_{\,2},\, b_{\,3},\, \cdots,\, b_{\,n}\,\right\| \,\leq\, 1\,,\;\; \left\|\,y,\, b_{\,2},\, b_{\,3},\, \cdots,\, b_{\,n}\,\right\| \,\leq\, 1\]
and 
\[\left\|\,x \,-\, y,\, b_{\,2},\, b_{\,3},\, \cdots,\, b_{\,n}\,\right\| \,\geq\, \epsilon.\]
Now, it is easily seen that
\[\left\|\,\dfrac{x_{\,k_{\,j}}}{1 \,+\, \dfrac{1}{j}},\, b_{\,2},\, b_{\,3},\, \cdots,\, b_{\,n}\,\right\| \,\leq\, 1\,,\; \,\left\|\,\dfrac{x_{\,k_{\,j \,+\, 1}}}{1 \,+\, \dfrac{1}{j}},\, b_{\,2},\, b_{\,3},\, \cdots,\, b_{\,n}\,\right\| \,\leq\, 1\]and
\[\left\|\,\dfrac{x_{\,k_{\,j \,+\, 1}} \,-\, x_{\,k_{\,j}}}{1 \,+\, \dfrac{1}{j}},\, b_{\,2},\, b_{\,3},\, \cdots,\, b_{\,n}\,\right\| \,\geq\, \dfrac{2\,\epsilon}{1 \,+\, \dfrac{1}{j}} \,\geq\, \epsilon.\]
Therefore
\[\left\|\,\dfrac{x_{\,k_{\,j \,+\, 1}} \,+\, x_{\,k_{\,j}}}{2\,\left(\,1 \,+\, \dfrac{1}{j}\,\right)},\, b_{\,2},\, b_{\,3},\, \cdots,\, b_{\,n}\,\right\| \,\leq\, 1 \,-\, \delta\]and this implies that
\[\lim\limits_{j \,\to\, \infty}\,\left\|\,x_{\,k_{\,j \,+\, 1}} \,+\, x_{\,k_{\,j}},\, b_{\,2},\, b_{\,3},\, \cdots,\, b_{\,n}\,\right\| \,\leq\, 2\,\left(\,1 \,-\, \delta\,\right) \,<\, 2.\]
This contradicts the assumption that 
\[\left\|\,x_{\,k} \,+\, x_{\,l},\, b_{\,2},\, b_{\,3},\, \cdots,\, b_{\,n}\,\right\| \,\to\, 2\; \;\;\textit{as}\; \;k,\, l \,\to\, \infty.\]Therefore, \,$\left\{\,x_{\,k}\,\right\}$\, is a Cauchy sequence and this proves the theorem.
\end{proof}

\section{Preliminaries of Banach Algebra in \,$n$-Banach space}
\smallskip\hspace{.6 cm}

In this section, we shall see that if a linear space is equipped with an algebra, many interesting results are obtained.\,We give the formal definition of a Banach algebra in \,$n$-Banach space.

\begin{definition}
Let \,$X$\, be a complex algebra of \,$\text{dim}\,X \,\geq\, n$\, with the \,$n$-norm \,$\left \|\,\cdot \,,\, \cdots \,,\, \cdot \,\right \|$.\,Then \,$X$\, is said to be a \,$n$-normed algebra if the \,$n$-norm satisfies the following inequality, known as multiplicative inequality
\begin{align}
 \left\|\, x\,y,\, a_{\,2},\, \cdots,\, a_{\,n} \,\right\| \,\leq\, \left\|\, x,\, a_{\,2},\, \cdots,\, a_{\,n}\,\right\|\,\left\|\,y,\, a_{\,2},\, \cdots,\, a_{\,n}\,\right\|\label{2.eq2.1}
\end{align}
for all \,$x,\, y,\, a_{\,2},\, a_{\,3},\, \cdots,\, a_{\,n} \,\in\, X$\, and if \,$X$\, contains a unit element \,$e$\, i.\,e., \,$x\,e \,=\, e\,x\, \,=\, x$, \,$x \,\in\, X$, and \,$ \left\|\,e,\, a_{\,2},\, \cdots,\, a_{\,n} \,\right\| \,=\, 1$, for all \,$a_{\,2},\, a_{\,3},\, \cdots,\, a_{\,n} \,\in\, X$.\,\,If, in addition, \,$X$\, is a \,$n$-Banach space with respect to the  \,$n$-norm \,$\left \|\,\cdot \,,\, \cdots \,,\, \cdot \,\right \|$\, then \,$X$\, is called a \,$n$-Banach algebra.   
\end{definition}

In the definition of a \,$n$-Banach algebra, we do not assume that \,$x\,y \,=\, y\,x$\, for all \,$x,\,y \,\in\, X$.\,But if this happens then \,$X$\, is called a commutative \,$n$-Banach algebra.\,A \,$n$-Banach subalgebra of the \,$n$-Banach algebra \,$X$\, is a subset of \,$X$\, which is itself a \,$n$-Banach algebra with respect to the algebraic operations defined on \,$X$\, and with the same identity and with the same \,$n$-norm. 

\begin{theorem}
Let \,$\left\{\,x_{\,k}\,\right\}$\, and \,$\left\{\,y_{\,k}\,\right\}$\, be two sequences in \,$X$\, such that \,$x_{\,k} \,\to\, x$\, and \,$y_{\,k} \,\to\, y$\, as \,$k \,\to\, \infty$, for some \,$x,\, y \,\in\, X$.\,Then \,$x_{\,k}\,y_{\,k} \,\to\, x\,y$\, as \,$k \,\to\, \infty$. 
\end{theorem}

\begin{proof}
Let \,$x,\, y,\, a_{\,2},\, \cdots,\, a_{\,n}  \,\in\, X$.\,Then,
\begin{align*}
&\left\|\,x_{\,k}\,y_{\,k} \,-\, x\,y \,,\, a_{\,2} \,,\, \cdots \,,\, a_{\,n} \,\right\|\\
&\left\|\,\left(\,x_{\,k} \,-\, x\,\right)\,y_{\,k} \,+\, x\,\left(\,y_{\,k} \,-\, y\,\right) \,,\, a_{\,2} \,,\, \cdots \,,\, a_{\,n} \,\right\|\\
&\leq\,\left\|\,x_{\,k} \,-\, x \,,\, a_{\,2} \,,\, \cdots \,,\, a_{\,n} \,\right\|\,\left\|\,y_{\,k} \,,\, a_{\,2} \,,\, \cdots \,,\, a_{\,n} \,\right\|\,+\\
&\hspace{.7cm} \,+\, \left\|\,y_{\,k} \,-\, y \,,\, a_{\,2} \,,\, \cdots \,,\, a_{\,n} \,\right\|\,\left\|\,x \,,\, a_{\,2} \,,\, \cdots \,,\, a_{\,n} \,\right\|\\
&\to\, 0\; \;\text{as}\; \, k \,\to\, \infty.
\end{align*}
Thus, \,$x_{\,k}\,y_{\,k} \,\to\, x\,y$\, as \,$k \,\to\, \infty$.\,This means that the multiplication is a continuous operation in \,$X$.   
\end{proof}

\begin{example}
Let \,$X$\, denote the set of all polynomials, 
\[x\,(\,t\,) \,=\, c_{\,0} \,+\, c_{\,1}\,t \,+\, c_{\,2}\,t^{\,2} \,+\, \,\cdots\, \,+\, c_{\,m}\,t^{\,m},\; \left(\,\;\text{$c_{\,0},\,\cdots,\,c_{\,m}$\, are complex}\,\right) \]
where \,$m \,\geq\, n$\, is not a fixed positive integer.\;Then \,$X$\, is a complex linear space with respect to the addition of  polynomials and scalar multiplication of a polynomial.\;For each \,$i \,=\, 1,\, 2,\, \cdots,\, n$, let \,$x_{\,i}\,(\,t\,) \,=\, a^{\,i}_{\,0} \,+\, a^{\,i}_{\,1}\,t \,+\, a^{\,i}_{\,2}\,t^{\,2} \,+\, \,\cdots\, \,+\, a^{\,i}_{\,m}\,t^{\,m}$.\;Now, define 
\[\left\|\,x_{\,1},\, \cdots,\, x_{\,n}\,\right\|\]
\begin{equation}\label{eq2.1}
\,=\, \begin{cases}
\max\limits_{\,j}\,\left|\,a^{\,1}_{\,j}\,\right| \,\times\,\cdots\,\times\max\limits_{\,j}\,\left|\,a^{\,n}_{\,j}\right| & \text{if}\; x_{\,1},\, \cdots \,,\, x_{\,n} \;\text{are linearly independent,} \\ 0 & \text{if}\; x_{\,1},\, \cdots \,,\, x_{\,n} \;\text{are linearly dependent}. \end{cases}
\end{equation} 
Let 
\[y_{\,1}\,(\,t\,) \,=\, b^{\,1}_{\,0} \,+\, b^{\,1}_{\,1}\,t \,+\, b^{\,1}_{\,2}\,t^{\,2} \,+\, \,\cdots\, \,+\, b^{\,1}_{\,m}\,t^{\,m} \,\in\, X\] and \,$\alpha \,\in\, \mathbb{K}$.\;Now, we see that
\begin{align*}
&\left\|\,\alpha\,x_{\,1},\, \cdots,\, x_{\,n}\,\right\| \,=\, \max\limits_{\,j}\,\left|\,\alpha\,a^{\,1}_{\,j}\,\right|\times\,\cdots\,\times\,\max\limits_{\,j}\,\left|\,a^{\,n}_{\,j}\right| \,=\, |\,\alpha\,|\,\left\|\,x_{\,1},\, \cdots,\, x_{\,n}\,\right\|,\\
& \left\|\,x_{\,1} \,+\, y_{\,1},\, \cdots,\, x_{\,n}\,\right\| \,=\, \max\limits_{\,j}\,\left|\,a^{\,1}_{\,j} \,+\, b^{\,1}_{\,j}\,\right|\times\,\cdots\times\,\max\limits_{\,j}\,\left|\,a^{\,n}_{\,j}\right|\\
&\leq\, \left(\,\max\limits_{\,j}\,\left|\,a^{\,1}_{\,j}\,\right| \,\times\,\cdots\,\times\max\limits_{\,j}\,\left|\,a^{\,n}_{\,j}\right|\,\right) \,+\, \left(\,\max\limits_{\,j}\,\left|\,b^{\,1}_{\,j}\,\right| \,\times\,\cdots\,\times\max\limits_{\,j}\,\left|\,a^{\,n}_{\,j}\right|\,\right)\\
&\,=\, \left\|\,x_{\,1},\, \cdots,\, x_{\,n}\,\right\| \,+\, \left\|\,y_{\,1},\, \cdots,\, x_{\,n}\,\right\|.
\end{align*}
Therefore, \,$X$\, becomes a linear\;$n$-normed space with respect to the \,$n$-norm defined by (\ref{eq2.1}).\,Now, we take the product of 
\begin{align*}
&x_{\,1}\,(\,t\,) \,=\, a^{\,1}_{\,0} \,+\, a^{\,1}_{\,1}\,t \,+\, a^{\,1}_{\,2}\,t^{\,2} \,+\, \,\cdots\, \,+\, a^{\,1}_{\,m}\,t^{\,m}\,\\ 
&y_{\,1}\,(\,t\,) \,=\, b^{\,1}_{\,0} \,+\, b^{\,1}_{\,1}\,t \,+\, b^{\,1}_{\,2}\,t^{\,2} \,+\, \,\cdots\, \,+\, b^{\,1}_{\,m}\,t^{\,m}
\end{align*}
as 
\[\left(\,x_{\,1}\,y_{\,1}\,\right)\,(\,t\,) \,=\, c^{\,1}_{\,0} \,+\, c^{\,1}_{\,1}\,t \,+\, c^{\,1}_{\,2}\,t^{\,2} \,+\, \,\cdots\, \,+\, c^{\,1}_{\,k}\,t^{\,k},\]
 where
 \[c_{\,k}^{\,1} \,=\, a^{\,1}_{\,0}\,b^{\,1}_{\,k} \,+\, a^{\,1}_{\,1}\,b^{\,1}_{\,k \,-\, 1} \,+\, \cdots \,+\, a^{\,1}_{\,k}\,b^{\,1}_{\,0}.\]Also, we have
\begin{align*}
& \left\|\,x_{\,1}\,y_{\,1},\, \cdots,\, x_{\,n}\,\right\| \,=\, \max\limits_{\,j}\,\left|\,a^{\,1}_{\,j}\;b^{\,1}_{\,j}\,\right|\times\,\cdots\times\,\max\limits_{\,j}\,\left|\,a^{\,n}_{\,j}\right|\\
&\leq\, \left(\,\max\limits_{\,j}\,\left|\,a^{\,1}_{\,j}\,\right| \,\times\,\cdots\,\times\max\limits_{\,j}\,\left|\,a^{\,n}_{\,j}\right|\,\right)\left(\,\max\limits_{\,j}\,\left|\,b^{\,1}_{\,j}\,\right| \,\times\,\cdots\,\times\max\limits_{\,j}\,\left|\,a^{\,n}_{\,j}\right|\,\right)\\
&\,=\, \left\|\,x_{\,1},\, \cdots,\, x_{\,n}\,\right\|\,\left\|\,y_{\,1},\, \cdots,\, x_{\,n}\,\right\|.
\end{align*}
Thus, \,$X$\, is commutative \,$n$-normed algebra with identity \,$(\,e \,=\, 1\,)$.
\end{example}

\begin{example}
The space \,$C\,[\,a,\, b\,]$\, with respect to the \,$n$-norm given by (\ref{eqpq2.1}) is a commutative \,$n$-Banach algebra with identity \,$e \,=\, 1$, where the product \,$x\,y$\, being defined as \,$(\,x\,y\,)\,(\,t\,) \,=\, x\,(\,t\,)\,y\,(\,t\,)$, for all \,$t \,\in\, [\,a,\, b\,]$.
\end{example}

\begin{definition}
Let \,$X$\, be a linear \,$n$-normed space and \,$b_{\,2},\,\cdots,\,b_{\,n} \,\in\, X$.\,Then an operator \,$T \,:\, X \,\to\, X$\, is called $b$-bounded if there exists some positive constant \,$M$\, such that
\[\left\|\,T\,x,\, b_{\,2},\,\cdots,\,b_{\,n}\,\right\| \,\leq\, M\,\left\|\,x,\, b_{\,2},\, \cdots,\, b_{\,n}\,\right\|\; \;\forall\; x \,\in\, X.\] 
\end{definition}

The norm of \,$T$\, is denoted by \,$\|\,T\,\|$\, and is defined as
\begin{align*}
&\|\,T\,\| \,=\, \inf\,\left\{\,M \,:\, \left\|\,T\,x,\, b_{\,2},\,\cdots,\,b_{\,n}\,\right\| \,\leq\, M\,\left\|\,x,\, b_{\,2},\, \cdots,\, b_{\,n}\,\right\|\; \;\forall\; x \,\in\, X\,\right\}.
\end{align*}

\begin{remark}
If \,$T$\, be a \,$b$-bounded on \,$X$, norm of \,$T$\, can be expressed by any one of the following equivalent formula:
\begin{itemize}
\item[(I)]\hspace{.2cm}$\|\,T\,\| \,=\, \sup\,\left\{\,\left\|\,T\,x,\, b_{\,2},\,\cdots,\,b_{\,n}\,\right\| \;:\; \left\|\,x,\, b_{\,2},\, \cdots,\, b_{\,n}\,\right\| \,\leq\, 1\,\right\}$.
\item[(II)]\hspace{.2cm}$\|\,T\,\| \,=\, \sup\,\left\{\,\left\|\,T\,x,\, b_{\,2},\,\cdots,\,b_{\,n}\,\right\| \;:\; \left\|\,x,\, b_{\,2},\, \cdots,\, b_{\,n}\,\right\| \,=\, 1\,\right\}$.
\item[(III)]\hspace{.2cm}$ \|\,T\,\| \,=\, \sup\,\left \{\,\dfrac{\left\|\,T\,x,\, b_{\,2},\,\cdots,\,b_{\,n}\,\right\|}{\left\|\,x,\, b_{\,2},\, \cdots,\, b_{\,n}\,\right\|} \;:\; \left\|\,x,\, b_{\,2},\, \cdots,\, b_{\,n}\,\right\| \,\neq\, 0\,\right \}$. 
\end{itemize}
Also, we have 
\[\left\|\,T\,x,\, b_{\,2},\,\cdots,\,b_{\,n}\,\right\| \,\leq\, \|\,T\,\|\, \left\|\,x,\, b_{\,2},\, \cdots,\, b_{\,n}\,\right\|\, \;\forall\; x \,\in\, X.\]
\end{remark}

\begin{example}
Let \,$X_{\,b}$\, be the set of all \,$b$-bounded linear operator on \,$H$, where \,$H$\, is a \,$n$-Hilbert space.\,We define the operations
\begin{align*}
&T_{\,1} \,+\, T_{\,2}\; \;\text{by}\; \;\left(\,T_{\,1} \,+\, T_{\,2}\,\right)\,(\,x\,) \,=\, T_{\,1}\,(\,x\,) \,+\, T_{\,2}\,(\,x\,)\,,\\
&\alpha\,T\; \;\text{by}\; \;\left(\,\alpha\,T\,\right)\,(\,x\,) \,=\, \alpha\,T\,(\,x\,)\,,\; \;T_{\,1}\,T_{\,2}\; \;\text{by}\; \;\left(\,T_{\,1}\,T_{\,2}\,\right)\,(\,x\,) \,=\, T_{\,1}\,\left(\,T_{\,2}\,(\,x\,)\,\right),
\end{align*}
for every \,$x \,\in\, H$.\,Define
\[\|\,T\,\| \,=\, \sup\,\left\{\,\left\|\,T\,x,\, b_{\,2},\,\cdots,\,b_{\,n}\,\right\| \;:\; \left\|\,x,\, b_{\,2},\, \cdots,\, b_{\,n}\,\right\| \,\leq\, 1\,\right\}.\]
Now, it can be easily verified that \,$\left(\,X_{\,b},\, \|\,\cdot\,\|\,\right)$\, is a Banach space.\,Also, we get that \,$\left\|\,T_{\,1}\,T_{\,2}\,\right\| \,\leq\, \left\|\,T_{\,1}\,\right\|\,\left\|\,T_{\,2}\,\right\|$.\,Thus, \,$X_{\,b}$\, is a Banach algebra.\,The identity operator is the identity element of this Banach algebra.\,Clearly, \,$X_{\,b}$\, is not commutative.    
\end{example}

If a \,$n$-Banach algebra \,$X$\, does not have an identity, we can extend \,$X$,\, which we show now, so that the extended \,$n$-Banach algebra has an identity.\,Let \,$\overline{X}$\, be the set of all ordered pairs \,$(\,x,\,a\,)$, where \,$x \,\in\, X$\, and \,$a \,\in\, \mathbb{C}$.\,Now, we show that \,$\overline{X}$\, can be equipped conveniently with algebraic operations and \,$n$-norm such that \,$\overline{X}$\, becomes a \,$n$-Banach algebra.\,We define the operations in \,$\overline{X}$\, as follows:
\begin{align*}
&(\,x,\,a\,) \,+\, (\,y,\,b\,) \,=\, \left(\,x \,+\, y,\, a \,+\, b\,\right)\,,\\
& \;(\,x,\,a\,)\,(\,y,\,b\,) \,=\, \left(\,x\,y \,+\, a\,y \,+\, b\,x,\, ab\,\right),\\ 
&\alpha\,(\,x,\,a\,) \,=\, (\,\alpha\,x,\,a\,\alpha\,)\,,\; \;\alpha \,\in\, \mathbb{K}.
\end{align*}   
We can easily verify that with the above operations, \,$\overline{X}$\, becomes an algebra.\,The \,$n$-norm on \,$\overline{X}$\, is defined as:
\[\left\|\,\left(\,x_{\,1},\, c_{\,1}\,\right),\, \left(\,x_{\,2},\, c_{\,2}\,\right),\, \cdots,\, \left(\,x_{\,n},\, c_{\,n}\,\right)\,\right\|\]
\[
\,=\, \begin{cases}
\left\|\, x_{\,1},\, x_{\,2},\, \cdots,\, x_{\,n}\,\right\| \,+\, \left|\,c_{\,1}\,c_{\,2}\, \cdots\, c_{\,n}\,\right| & \text{if}\; (\,x_{\,1},\, c_{\,1}\,),\, \cdots,\, (\,x_{\,n},\, c_{\,n}\,) \;\text{are linearly,}\\ &\hspace{1.5cm}\text{independent} \\ 0 & \text{if}\; (\,x_{\,1},\, c_{\,1}\,),\, \cdots,\, (\,x_{\,n},\, c_{\,n}\,) \;\text{are linearly}\\ &\hspace{1.5cm}\text{dependent}, \end{cases}\]
where \,$x_{\,1},\, x_{\,2},\, \cdots,\, x_{\,n} \,\in\, X$\, and \,$c_{\,1},\, c_{\,2},\, \cdots,\, c_{\,n} \,\in\, \mathbb{C}$.\,Here, the \,$n$-norm axioms are fulfilled.\,Also, for \,$x,\, y,\, a_{\,2},\, \cdots,\, a_{\,n} \,\in\, X$\, and \,$a,\, b,\, b_{\,2},\, \cdots,\, b_{\,n} \,\in\, \mathbb{C}$, we have 
\begin{align*}
&\left\|\,\left(\,x,\, a\,\right)\,\left(\,y,\, b\,\right),\, \left(\,a_{\,2},\, b_{\,2}\,\right),\, \cdots,\, \left(\,a_{\,n},\, b_{\,n}\,\right)\,\right\| \\
&\,=\, \left\|\,\left(\,x\,y \,+\, a\,y \,+\, b\,x,\, ab\,\right),\, \left(\,a_{\,2},\, b_{\,2}\,\right),\, \cdots,\, \left(\,a_{\,n},\, b_{\,n}\,\right)\,\right\|\\
&=\, \left\|\,x\,y \,+\, a\,y \,+\, b\,x,\, a_{\,2},\, \cdots,\, a_{\,n}\,\right\| \,+\, \left|\,a\,b\,b_{\,2}\,\cdots\,b_{\,n}\,\right|\\
&\leq\,\left\|\,x,\, a_{\,2},\, \cdots,\, a_{\,n}\,\right\|\,\left\|\,y,\, a_{\,2},\, \cdots,\, a_{\,n}\,\right\| \,+\, |\,a\,|\left\|\,y,\, a_{\,2},\, \cdots,\, a_{\,n}\,\right\|\,+\\
&\hspace{1cm} \,+\, |\,b\,|\left\|\,x,\, a_{\,2},\, \cdots,\, a_{\,n}\,\right\| \,+\, \left|\,a\,b\,b_{\,2}\,\cdots\,b_{\,n}\,\right|\\
&\leq\,\left(\,\left\|\,x,\, a_{\,2},\, \cdots,\, a_{\,n}\,\right\| \,+\, \left|\,a\,b_{\,2}\,\cdots\,b_{\,n}\,\right|\,\right)\,\left(\,\left\|\,y,\, a_{\,2},\, \cdots,\, a_{\,n}\,\right\| \,+\, \left|\,b\,b_{\,2}\,\cdots\,b_{\,n}\,\right|\,\right)\\
&=\,\left\|\,\left(\,x,\, a\,\right),\, \left(\,a_{\,2},\, b_{\,2}\,\right),\, \cdots,\, \left(\,a_{\,n},\, b_{\,n}\,\right)\,\right\|\,\left\|\,\left(\,y,\, b\,\right),\, \left(\,a_{\,2},\, b_{\,2}\,\right),\, \cdots,\, \left(\,a_{\,n},\, b_{\,n}\,\right)\,\right\|.
\end{align*}
We now show that \,$\overline{\,X}$\, is a \,$n$-Banach space i.\,e., \,$\overline{\,X}$\, is complete.\,Clearly, a sequence \,$\left\{\,\left(\,x_{\,k},\, a_{\,k}\,\right)\,\right\}$\, is a Cauchy sequence in \,$\overline{\,X}$\, if and only if \,$\left\{\,x_{\,k}\,\right\}$\, is a Cauchy sequence in \,$X$\, and \,$\left\{\,a_{\,k}\,\right\}$\, is a Cauchy sequence in \,$\mathbb{C}$.\,Since \,$X$\, and \,$\mathbb{C}$\, are complete, there exist \,$x \,\in\, X$\, and \,$a \,\in\, \mathbb{C}$\, such that \,$\lim\limits_{k \,\to\, \infty}\,x_{\,k} \,=\, x$\, and \,$\lim\limits_{k \,\to\, \infty}\,a_{\,k} \,=\, a$.\,Therefore, \,$\lim\limits_{k \,\to\, \infty}\,\left(\,x_{\,k},\, a_{\,k}\,\right) \,=\, (\,x,\, a\,)$\, is in \,$\overline{\,X}$.\,Also, it is easy to verify that \,$\left(\,x,\, a\,\right)\,\left(\,\theta,\, 1\,\right) \,=\, \left(\,x,\, a\,\right) \,=\, \left(\,\theta,\, 1\,\right)\,\left(\,x,\, a\,\right)$, where \,$\theta$\, is the zero element in \,$X$.\,Thus, \,$\overline{\,X}$\, is a \,$n$-Banach algebra with the indentity \,$\left(\,\theta,\, 1\,\right)$.   

For $2$-normed algebra augmentation of unity is discussed by N. Srivastava et al. \cite{NSS}.

\subsection{Invertible and Non-invertible Elements}
\smallskip\hspace{.6 cm}

In this section, we deal with the set of invertible and non-invertible elements in a \,$n$-Banach algebra.\,We remember that \,$G$\, denotes the set of all invertible elements of \,$X$.\,It is easy to verify that \,$G$\, is a group.\,Let \,$S$\, denotes the set of all non-invertible elements of \,$X$.\,Clearly, \,$S \,\cup\, G \,=\, X$.

\begin{theorem}\label{3.th3.11}
Let \,$X$\, be a complex \,$n$-Banach algebra with unity \,$e$.\,If \,$x \,\in\, X$\, satisfies \,$\left\|\,x,\, a_{\,2},\, \cdots,\, a_{\,n}\,\right\| \,<\, 1$, for every \,$a_{\,2},\, \cdots,\, a_{\,n}  \,\in\, X$, then \,$e \,-\, x$\, is invertible and \,$\left(\,e \,-\, x\,\right)^{\,-\, 1} \,=\, e \,+\, \sum\limits_{j \,=\, 1}^{\,\infty}\,x^{\,j}$.
\end{theorem}

\begin{proof}
From (\ref{2.eq2.1}), we get
\[\left\|\,x^{\,j},\, a_{\,2},\, \cdots,\, a_{\,n}\,\right\| \,\leq\, \left\|\,x,\, a_{\,2},\, \cdots,\, a_{\,n}\,\right\|^{\,j},\;\; \;\text{for every}\; \;a_{\,2},\, \cdots,\, a_{\,n}  \,\in\, X,\]
for any positive integer \,$j$, so that the infinite series \,$\sum\limits_{j \,=\, 1}^{\,\infty}\,x^{\,j}$\, is convergent because \,$\left\|\,x,\, a_{\,2},\, \cdots,\, a_{\,n}\,\right\| \,<\, 1$.\,Since \,$X$\, is complete, it is easy to verify that the infinite series \,$\sum\limits_{j \,=\, 1}^{\,\infty}\,x^{\,j}$\, is convergent to some element in \,$X$.\,Let therefore \,$s \,=\, e \,+\, \sum\limits_{j \,=\, 1}^{\,\infty}\,x^{\,j}$.\,If we prove that \,$s \,=\, \left(\,e \,-\, x\,\right)^{\,-\, 1}$, the proof is finished.\,Now, by a simple calculation, we obtain
\begin{align}
\left(\,e \,-\, x\,\right)\,\left(\,e \,+\, x \,+\, \cdots \,+\, x^{\,k}\,\right)& \,=\, \left(\,e \,+\, x \,+\, \cdots \,+\, x^{\,k}\,\right)\,\left(\,e \,-\, x\,\right)\nonumber\\
& \,=\, e \,-\, x^{\,k \,+\, 1}.\label{3.eq3.1} 
\end{align} 
Since \,$\left\|\,x,\, a_{\,2},\, \cdots,\, a_{\,n}\,\right\| \,<\, 1$\, and so \,$x^{\,k \,+\, 1} \,\to\, \theta$\, as \,$k \,\to\, \infty$.\,Therefore, letting \,$k \,\to\, \infty$\, in (\ref{3.eq3.1}) and using the continuity of multiplication in \,$X$, we obtain \,$\left(\,e \,-\, x\,\right)\,s \,=\, s\,\left(\,e \,-\, x\,\right) \,=\, e$.\,Thus, \,$s \,=\, \left(\,e \,-\, x\,\right)^{\,-\, 1}$.\,This proves the theorem.     
\end{proof}

In Theorem \ref{3.th3.11}, if \,$x$\, is replaced by \,$\left(\,e \,-\, x\,\right)$, we obtain the following corollary.

\begin{corollary}\label{3.cor3.1}
If \,$x \,\in\, X$\, satisfies \,$\left\|\,e \,-\, x,\, a_{\,2},\, \cdots,\, a_{\,n}\,\right\| \,<\, 1$, for every \,$a_{\,2},\, \cdots,\\ a_{\,n}  \,\in\, X$, then \,$x^{\,-\, 1}$\, exists and \,$x^{\,-\, 1} \,=\, e \,+\, \sum\limits_{j \,=\, 1}^{\,\infty}\,\left(\,e \,-\, x\,\right)^{\,j}$.
\end{corollary}

Next, we prove another corollary.

\begin{corollary}\label{3.cor3.2}
Let \,$x \,\in\, X$\, and \,$\lambda$\, be a scalar such that \,$\left\|\,x,\, a_{\,2},\, \cdots,\, a_{\,n}\,\right\| \,<\, |\,\lambda\,|$, for every \,$a_{\,2},\, \cdots,\, a_{\,n} \,\in\, X$.\,Then \,$\left(\,\lambda\,e \,-\, x\,\right)^{\,-\, 1}$\, exists and \[\left(\,\lambda\,e \,-\, x\,\right)^{\,-\, 1} \,=\, \sum\limits_{j \,=\, 1}^{\,\infty}\,\lambda^{\,-\, j}\,x^{\,j \,-\, 1}, \;\left(\,x^{\,0} \,=\, e\,\right)\]
\end{corollary}

\begin{proof}
If \,$y \,\in\, X$\, be such that \,$y^{\,-\, 1} \,\in\, X$\, and \,$\alpha \,\neq\, 0$\, be a scalar, then it is clear that \,$\left(\,\alpha\,y\,\right)^{\,-\, 1}$\, exists and \,$\left(\,\alpha\,y\,\right)^{\,-\, 1} \,=\, \alpha^{\,-\, 1}\,y^{\,-\, 1}$.\,Having noted this, we can write \,$\lambda\,e \,-\, x \,=\, \lambda\,\left(\,e \,-\, \dfrac{x}{\lambda}\,\right)$\, and we show that \,$\left(\,e \,-\, \dfrac{x}{\lambda}\,\right)^{\,-\, 1}$\, exists.\,Now, by hypothesis, we have
\begin{align*}
\left\|\,e \,-\, \left(\,e \,-\, \dfrac{x}{\lambda}\,\right),\, a_{\,2},\, \cdots,\, a_{\,n}\,\right\| \,=\, \dfrac{\left\|\,x,\, a_{\,2},\, \cdots,\, a_{\,n}\,\right\|}{\left|\,\lambda\,\right|} \,<\, 1.
\end{align*} 
So, by Corollary \ref{3.cor3.1}, \,$\left(\,e \,-\, \dfrac{x}{\lambda}\,\right)^{\,-\, 1}$\, exists and so \,$\left(\,\lambda\,e \,-\, x\,\right)^{\,-\, 1}$\, exists.

For the infinite series representation, we make use of the second part of Corollary \ref{3.cor3.1}.\,Now, we have \,$\left(\,\lambda\,e \,-\, x\,\right) \,=\, \lambda\,\left(\,e \,-\, \lambda^{\,-\, 1}\,x\,\right)$\, and it is easy to verify that \,$\left(\,\lambda\,e \,-\, x\,\right)^{\,-\, 1} \,=\, \sum\limits_{j \,=\, 1}^{\,\infty}\,\lambda^{\,-\, j}\,x^{\,j \,-\, 1}$.\,This proves the Corollary.     
\end{proof}

In the following Theorem, we obtain a set-theoretic property of invertible and non-invertible elements of \,$X$.

\begin{theorem}
The set \,$G$\, of all invertible elements of \,$X$\, is open subset.
\end{theorem}

\begin{proof}
Let \,$x_{\,0} \,\in\, G$.\,We should show that there exists a sphere around \,$x_{\,0}$, each element of which belongs to \,$G$.\,Consider the open sphere \[B_{\,\left\{\,a_{\,2} \,,\, \cdots \,,\, a_{\,n}\,\right\}}\,\left(\,x_{\,0} \,,\, \dfrac{1}{\left\|\,x^{\,-\, 1}_{\,0},\, a_{\,2},\, \cdots,\, a_{\,n}\,\right\|}\,\right)\]
with center at \,$x_{\,0}$\, and radius \,$1 \,/\, \left\|\,x^{\,-\, 1}_{\,0},\, a_{\,2},\, \cdots,\, a_{\,n}\,\right\|$.\,For every point \,$x$\, of this sphere satisfies the inequality 
\begin{align}
\left\|\,x \,-\, x_{\,0},\, a_{\,2},\, \cdots,\, a_{\,n}\,\right\| \,<\, \dfrac{1}{\left\|\,x^{\,-\, 1}_{\,0},\, a_{\,2},\, \cdots,\, a_{\,n}\,\right\|}\label{3.eq3.2} 
\end{align}
Let \,$y \,=\, x^{\,-\, 1}_{\,0}\,x$\, and \,$z \,=\, e \,-\, y$.\,Then by (\ref{2.eq2.1}) and (\ref{3.eq3.2}), we have
\begin{align*}
\left\|\,z,\, a_{\,2},\, \cdots,\, a_{\,n}\,\right\|& \,=\, \left\|\,y \,-\, e,\, a_{\,2},\, \cdots,\, a_{\,n}\,\right\| \,=\, \left\|\,x^{\,-\, 1}_{\,0}\,x \,-\, x^{\,-\, 1}_{\,0}\,x_{\,0},\, a_{\,2},\, \cdots,\, a_{\,n}\,\right\|\\
&\,=\,\left\|\,x^{\,-\, 1}_{\,0}\,\left(\,x \,-\, x_{\,0}\,\right),\, a_{\,2},\, \cdots,\, a_{\,n}\,\right\|\\
&\leq\,\left\|\,x^{\,-\, 1}_{\,0},\, a_{\,2},\, \cdots,\, a_{\,n}\,\right\|\,\left\|\,x \,-\, x_{\,0},\, a_{\,2},\, \cdots,\, a_{\,n}\,\right\| \,<\, 1.
\end{align*}
By Theorem \ref{3.th3.11}, \,$e \,-\, z$\, is invertible i.\,e., \,$y$\, is invertible.\,So, \,$y \,\in\, G$.\,Now, \,$x_{\,0} \,\in\, G\,,\; \;y \,\in\, G$\, and by our earlier verification, \,$G$\, is group, so \,$x_{\,0}\,y \,\in\, G$.\,But \,$x_{\,0}\,y \,=\, x_{\,0}\,x^{\,-\, 1}\,x \,=\, x$.\,So, any \,$x$\, satisfying (\ref{3.eq3.2}) belongs to \,$G$.\,This shows that \,$G$\, is an open subset of \,$X$.\,This completes the proof.    
\end{proof}

\begin{corollary}\label{3.cor3.22}
The set of all non-invertible elements is a closed subset of \,$X$.
\end{corollary}

\begin{definition}
A linear operator \,$T \,:\, X \,\to\, X$\, is said to be continuous at \,$x_{\,0} \,\in\, X$\, if for any open ball \,$B_{\,\{\,e_{\,2},\, \cdots,\, e_{\,n}\,\}}\,(\,T\,(\,x_{\,0}\,),\, \epsilon\,) \,>\, 0$\, there exists a \,$\delta \,>\, 0$\, such that
\[T\,\left(\,B_{\,\{\,e_{\,2},\, \cdots,\, e_{\,n}\,\}}\,(\,x_{\,0},\, \delta\,)\,\right) \,\subset\, B_{\,\{\,e_{\,2},\, \cdots,\, e_{\,n}\,\}}\,(\,T\,(\,x_{\,0}\,),\, \epsilon\,).\]
Equivalently, for given \,$\epsilon \,>\, 0$, there exist some \,$e_{\,2},\, \cdots,\, e_{\,n} \,\in\, X$\, and \,$\delta \,>\, 0$\, such that for \,$x \,\in\, X$
\[\left\|\,x \,-\, x_{\,0},\,e_{\,2},\, \cdots,\, e_{\,n}\,\right\| \,<\, \delta \,\Rightarrow\, \left\|\,T\,(\,x\,) \,-\, T\,(\,x_{\,0}\,),\,e_{\,2},\, \cdots,\, e_{\,n}\,\right\| \,<\, \epsilon.\]
\end{definition}

\begin{theorem}
The mapping \,$T \,:\, G \,\to\, G$\, given by \,$T\,(\,x\,) \,=\, x^{\,-\, 1}$, for all \,$x \,\in\, G$\, is continuous.
\end{theorem}

\begin{proof}
Let \,$x_{\,0}$\, be an element of \,$G$\, and \,$x$\, be any other element of \,$G$\, that satisfies the inequality
\begin{align*}
\left\|\,x \,-\, x_{\,0},\, a_{\,2},\, \cdots,\, a_{\,n}\,\right\| \,<\, \dfrac{1}{2\,\left\|\,x^{\,-\, 1}_{\,0},\, a_{\,2},\, \cdots,\, a_{\,n}\,\right\|}.
\end{align*}
Now,
\begin{align}
&\left\|\,x^{\,-\, 1}_{\,0}\,x \,-\, e,\, a_{\,2},\, \cdots,\, a_{\,n}\,\right\| \,=\, \left\|\,x^{\,-\, 1}_{\,0}\,\left(\,x \,-\, x_{\,0}\,\right),\, a_{\,2},\, \cdots,\, a_{\,n}\,\right\|\nonumber\\ 
&\,\leq\,\left\|\,x^{\,-\, 1}_{\,0},\, a_{\,2},\, \cdots,\, a_{\,n}\,\right\|\,\left\|\,x \,-\, x_{\,0},\, a_{\,2},\, \cdots,\, a_{\,n}\,\right\| \,<\, \dfrac{1}{2} \label{3.eq3.34}
\end{align}
and so by Corollary \ref{3.cor3.1}, \,$x^{\,-\, 1}_{\,0}\,x \,\in\, G$\, and further
\begin{align}
x^{\,-\, 1}\,x_{\,0} \,=\, \left(\,x^{\,-\, 1}_{\,0}\,x\,\right)^{\,-\, 1} \,=\, e \,+\, \sum\limits_{j \,=\, 1}^{\,\infty}\,\left(\,e \,-\, x^{\,-\, 1}_{\,0}\,x\,\right)^{\,j}.\label{3.eq3.35}
\end{align}
Next, we show that the mapping \,$T \,:\, G \,\to\, G$\, given by \,$T\,(\,x\,) \,=\, x^{\,-\, 1}$\, is continuous at \,$x_{\,0}$.\,Now, we have
\begin{align}
&\left\|\,T\,x \,-\, T\,x_{\,0},\, a_{\,2},\, \cdots,\, a_{\,n}\,\right\| \,=\, \left\|\,x^{\,-\, 1} \,-\, x^{\,-\, 1}_{\,0},\, a_{\,2},\, \cdots,\, a_{\,n}\,\right\|\nonumber\\
&=\,\left\|\,\left(\,x^{\,-\, 1}\,x_{\,0} \,-\, e\,\right)\,x^{\,-\, 1}_{\,0},\, a_{\,2},\, \cdots,\, a_{\,n}\,\right\|\nonumber\\
&\leq\,\left\|\,x^{\,-\, 1}_{\,0},\, a_{\,2},\, \cdots,\, a_{\,n}\,\right\|\,\left\|\,x^{\,-\, 1}\,x_{\,0} \,-\, e,\, a_{\,2},\, \cdots,\, a_{\,n}\,\right\|\nonumber\\
&=\,\left\|\,x^{\,-\, 1}_{\,0},\, a_{\,2},\, \cdots,\, a_{\,n}\,\right\|\,\left\|\,\sum\limits_{j \,=\, 1}^{\,\infty}\,\left(\,e \,-\, x^{\,-\, 1}_{\,0}\,x\,\right)^{\,j},\, a_{\,2},\, \cdots,\, a_{\,n}\,\right\|\;[\,\text{by}\; (\ref{3.eq3.35})\;]\nonumber\\
&\leq\,\left\|\,x^{\,-\, 1}_{\,0},\, a_{\,2},\, \cdots,\, a_{\,n}\,\right\|\,\sum\limits_{j \,=\, 1}^{\,\infty}\,\left\|\,e \,-\, x^{\,-\, 1}_{\,0}\,x,\, a_{\,2},\, \cdots,\, a_{\,n}\,\right\|^{\,j}\nonumber\\
&=\,\left\|\,x^{\,-\, 1}_{\,0},\, a_{\,2},\, \cdots,\, a_{\,n}\,\right\|\,\left\|\,e \,-\, x^{\,-\, 1}_{\,0}\,x,\, a_{\,2},\, \cdots,\, a_{\,n}\,\right\|\sum\limits_{j \,=\, 0}^{\,\infty}\,\left\|\,e \,-\, x^{\,-\, 1}_{\,0}\,x,\, a_{\,2},\, \cdots,\, a_{\,n}\,\right\|^{\,j}\nonumber\\
&=\,\dfrac{\left\|\,x^{\,-\, 1}_{\,0},\, a_{\,2},\, \cdots,\, a_{\,n}\,\right\|\,\left\|\,e \,-\, x^{\,-\, 1}_{\,0}\,x,\, a_{\,2},\, \cdots,\, a_{\,n}\,\right\|}{1 \,-\, \left\|\,e \,-\, x^{\,-\, 1}_{\,0}\,x,\, a_{\,2},\, \cdots,\, a_{\,n}\,\right\|}\nonumber\\
&<\, 2\,\left\|\,x^{\,-\, 1}_{\,0},\, a_{\,2},\, \cdots,\, a_{\,n}\,\right\|\,\left\|\,e \,-\, x^{\,-\, 1}_{\,0}\,x,\, a_{\,2},\, \cdots,\, a_{\,n}\,\right\|\; \;[\;\text{by}\; (\ref{3.eq3.34})\;]\nonumber\\
&\leq\,2\,\left\|\,x^{\,-\, 1}_{\,0},\, a_{\,2},\, \cdots,\, a_{\,n}\,\right\|^{\,2}\,\left\|\,x \,-\, x_{\,0},\, a_{\,2},\, \cdots,\, a_{\,n}\,\right\|,\label{3.eq3.36}
\end{align}
since 
\begin{align*}
&\left\|\,e \,-\, x^{\,-\, 1}_{\,0}\,x,\, a_{\,2},\, \cdots,\, a_{\,n}\,\right\| \,=\, \left\|\,x^{\,-\, 1}_{\,0}\,x \,-\, e,\, a_{\,2},\, \cdots,\, a_{\,n}\,\right\|\\
&=\,\left\|\,x^{\,-\, 1}_{\,0}\,\left(\,x \,-\, x_{\,0}\,\right)\, a_{\,2},\, \cdots,\, a_{\,n}\,\right\|\\
&\leq\,\left\|\,x^{\,-\, 1}_{\,0},\, a_{\,2},\, \cdots,\, a_{\,n}\,\right\|\,\left\|\,x \,-\, x_{\,0},\, a_{\,2},\, \cdots,\, a_{\,n}\,\right\|.
\end{align*}
Thus, from (\ref{3.eq3.36}), the continuity of \,$T$\, at \,$x_{\,0}$\, follows.\,This completes the proof. 
\end{proof}

\begin{theorem}
Let \,$x \,\in\, G$\, and \,$h \,\in\, X$\, be such that
\begin{align*}
\left\|\,h,\, a_{\,2},\, \cdots,\, a_{\,n}\,\right\| \,<\, \dfrac{1}{2}\,\left\|\,x^{\,-\, 1},\, a_{\,2},\, \cdots,\, a_{\,n}\,\right\|^{\,-\, 1},
\end{align*}
for every \,$a_{\,2},\, \cdots,\, a_{\,n} \,\in\, X$.\,Then \,$x \,+\, h \,\in\, G$\, and 
\begin{align*}
&\left\|\,\left(\,x \,+\, h\,\right)^{\,-\, 1} \,-\, x^{\,-\, 1} \,+\, x^{\,-\, 1}\,h\,x^{\,-\, 1},\, a_{\,2},\, \cdots,\, a_{\,n}\,\right\|\\
&\leq\,2\,\left\|\,x^{\,-\, 1},\, a_{\,2},\, \cdots,\, a_{\,n}\,\right\|^{\,3} \,\left\|\,h,\, a_{\,2},\, \cdots,\, a_{\,n}\,\right\|^{\,2}.  
\end{align*}
\end{theorem}

\begin{proof}
Since for every \,$a_{\,2},\, \cdots,\, a_{\,n} \,\in\, X$, we have 
\[\left\|\,h,\, a_{\,2},\, \cdots,\, a_{\,n}\,\right\|\,\left\|\,x^{\,-\, 1},\, a_{\,2},\, \cdots,\, a_{\,n}\,\right\| \,<\,  \dfrac{1}{2},\]
 and therefore
\begin{align*}
\left\|\,x^{\,-\, 1}\,h,\, a_{\,2},\, \cdots,\, a_{\,n}\,\right\| \,\leq\, \left\|\,x^{\,-\, 1},\, a_{\,2},\, \cdots,\, a_{\,n}\,\right\|\,\left\|\,h,\, a_{\,2},\, \cdots,\, a_{\,n}\,\right\| \,<\, \dfrac{1}{2}. 
\end{align*}
Thus, by Theorem \ref{3.th3.11}, we get that \,$e \,+\, x^{\,-\, 1}\,h \,\in\, G$.\,Since \,$G$\, is a group, we have that \,$x\,\left(\,e \,+\, x^{\,-\, 1}\,h\,\right) \,\in\, G$.\,But \,$x\,\left(\,e \,+\, x^{\,-\, 1}\,h\,\right) \,=\, x \,+\, h$\, and so \,$x \,+\, h \,\in\, G$.

Before proving the second part, we note by Theorem \ref{3.th3.11} that 
\[\left(\,e \,+\, x^{\,-\, 1}\,h\,\right)^{\,-\, 1} \,=\, e \,+\, \sum\limits_{j \,=\, 1}^{\,\infty}\,(\,-\, 1\,)^{\,j}\,\left(\,x^{\,-\, 1}\,h\,\right)^{\,j}\]
and so
\begin{align}
&\left\|\,\left(\,e \,+\, x^{\,-\, 1}\,h\,\right)^{\,-\, 1} \,-\, e \,+\, x^{\,-\, 1}\,h,\, a_{\,2},\, \cdots,\, a_{\,n}\,\right\|\nonumber\\
&\leq\,\sum\limits_{j \,=\, 2}^{\,\infty}\,\left\|\,x^{\,-\, 1}\,h,\, a_{\,2},\, \cdots,\, a_{\,n}\,\right\|^{\,j}\nonumber\\
&=\,\left\|\,x^{\,-\, 1}\,h,\, a_{\,2},\, \cdots,\, a_{\,n}\,\right\|^{\,2}\,\sum\limits_{j \,=\, 0}^{\,\infty}\,\left\|\,x^{\,-\, 1}\,h,\, a_{\,2},\, \cdots,\, a_{\,n}\,\right\|^{\,j}\nonumber\\
&<\,\left\|\,x^{\,-\, 1}\,h,\, a_{\,2},\, \cdots,\, a_{\,n}\,\right\|^{\,2}\,\sum\limits_{j \,=\, 0}^{\,\infty}\,\dfrac{1}{2^{\,j}} \,=\,2\,\left\|\,x^{\,-\, 1}\,h,\, a_{\,2},\, \cdots,\, a_{\,n}\,\right\|^{\,2}.\label{3.eq3.37}
\end{align}
Now,
\begin{align*}
&\left(\,x \,+\, h\,\right)^{\,-\, 1} \,-\, x^{\,-\, 1} \,+\, x^{\,-\, 1}\,h\,x^{\,-\, 1}\\
& \,=\, \left[\,x\,\left(\,e \,+\, x^{\,-\, 1}\,h\,\right)\,\right]^{\,-\, 1} \,-\, x^{\,-\, 1} \,+\, x^{\,-\, 1}\,h\,x^{\,-\, 1}\\
&=\,\left(\,e \,+\, x^{\,-\, 1}\,h\,\right)^{\,-\, 1}\,x^{\,-\, 1} \,-\, x^{\,-\, 1} \,+\, x^{\,-\, 1}\,h\,x^{\,-\, 1}\\
&=\,\left[\,\left(\,e \,+\, x^{\,-\, 1}\,h\,\right)^{\,-\, 1} \,-\, e \,+\, x^{\,-\, 1}\,h\,\right]\,x^{\,-\, 1}   
\end{align*}
and so taking \,$n$-norm and using (\ref{3.eq3.37}), for every \,$a_{\,2},\, \cdots,\, a_{\,n} \,\in\, X$, we get
\begin{align*}
&\left\|\,\left(\,x \,+\, h\,\right)^{\,-\, 1} \,-\, x^{\,-\, 1} \,+\, x^{\,-\, 1}\,h\,x^{\,-\, 1},\, a_{\,2},\, \cdots,\, a_{\,n}\,\right\|\\
&\leq\,\left\|\,\left(\,e \,+\, x^{\,-\, 1}\,h\,\right)^{\,-\, 1} \,-\, e \,+\, x^{\,-\, 1}\,h,\, a_{\,2},\, \cdots,\, a_{\,n}\,\right\|\,\left\|\,x^{\,-\, 1},\, a_{\,2},\, \cdots,\, a_{\,n}\,\right\|\\
&\leq\,2\,\left\|\,x^{\,-\, 1}\,h,\, a_{\,2},\, \cdots,\, a_{\,n}\,\right\|^{\,2}\,\left\|\,x^{\,-\, 1},\, a_{\,2},\, \cdots,\, a_{\,n}\,\right\|\\
&\leq\,2\,\left\|\,x^{\,-\, 1},\, a_{\,2},\, \cdots,\, a_{\,n}\,\right\|^{\,3}\,\left\|\,h,\, a_{\,2},\, \cdots,\, a_{\,n}\,\right\|.
\end{align*}
This proves the Theorem.
\end{proof}

\subsection{Topological Divisor of Zero}
\smallskip\hspace{.6 cm}

In this section, we give the idea of a topological divisor of zero in a \,$n$-Banach algebra and establish its relationship with non-inevrtible elements.

\begin{definition}
Let \,$X$\, be a \,$n$-Banach algebra.\,An element \,$z \,\in\, X$\, is said to be a topological divisor of zero if there exists a sequence \,$\left\{\,z_{\,k}\,\right\}$\, in \,$X$\, with \,$\left\|\,z_{\,k},\, a_{\,2},\, \cdots,\, a_{\,n}\,\right\| \,=\, 1$, for every \,$a_{\,2},\, \cdots,\, a_{\,n} \,\in\, X$\, and \,$k \,=\, 1,\, 2,\, \cdots$\, and such that either \,$z\,z_{\,k} \,\to\, \theta$\, or \,$z_{\,k}\,z \,\to\, \theta$\, as \,$k \,\to\, \infty$.
\end{definition}

Clearly, every divisor of zero is also a topological divisor of zero.\,Let \,$Z$\, denote the set of all topological divisors of zero in \,$X$.\,In the next Theorem, we gives certain connection between \,$Z$\, and the set \,$S$\, of all non-invertible elements.

\begin{theorem}
The set \,$Z$\, is a subset of \,$S$.
\end{theorem}    

\begin{proof}
Let \,$z \,\in\, Z$.\,Then there exists a sequence \,$\left\{\,z_{\,k}\,\right\}$\, with \,$\left\|\,z_{\,k},\, a_{\,2},\, \cdots,\, a_{\,n}\,\right\| \,=\, 1$, for every \,$a_{\,2},\, \cdots,\, a_{\,n} \,\in\, X$\, and \,$k \,=\, 1,\, 2,\, \cdots$\, and such that either \,$z\,z_{\,k} \,\to\, \theta$\, or \,$z_{\,k}\,z \,\to\, \theta$\, as \,$k \,\to\, \infty$.\,Suppose that \,$z\,z_{\,k} \,\to\, \theta$\, as \,$k \,\to\, \infty$.\,If possible, let \,$z \,\in\, G$.\,Then \,$z^{\,-\, 1}$\, exists.\,Now, as multiplication is a continuous operation, we should have \,$z_{\,k} \,=\, z^{\,-\, 1}\,\left(\,z\,z_{\,k}\,\right) \,\to\, z^{\,-\, 1}\,\theta \,=\, \theta$\, as \,$k \,\to\, \infty$.\,This contradicts the fact that \,$\left\|\,z_{\,k},\, a_{\,2},\, \cdots,\, a_{\,n}\,\right\| \,=\, 1$, for every \,$a_{\,2},\, \cdots,\, a_{\,n} \,\in\, X$\, and \,$k \,=\, 1,\, 2,\, \cdots$.\,So \,$z \,\in\, S$\, and this proves the Theorem.  
\end{proof}

\begin{theorem}
The boundary of \,$S$\, is a subset of \,$Z$.
\end{theorem}

\begin{proof}
By Corollary \ref{3.cor3.22}, \,$S$\, is a closed subset of \,$X$.\,So, any boundary point of \,$S$\, is in \,$S$.\,Further, for every boundary point of \,$S$, there exists a sequence of elements from \,$G$\, that converges to the boundary point.\,Let \,$x$\, be any boundary point of \,$S$, then \,$x \,\in\, S$\, and there exists a sequence \,$\left\{\,s_{\,k}\,\right\}$\, of elements from \,$G$\, such that \,$s_{\,k} \,\to\, x$\, as \,$k \,\to\, \infty$.\,Now, we have \,$s^{\,-\, 1}_{\,k}\,x \,-\, e \,=\, s^{\,-\, 1}_{\,k}\,\left(\,x \,-\, s_{\,k}\,\right)$.\,Therefore, if \,$\left\|\,s^{\,-\, 1}_{\,k},\, a_{\,2},\, \cdots,\, a_{\,n}\,\right\|$\, is a bounded sequence then because \,$s_{\,k} \,\to\, x$\, as \,$k \,\to\, \infty$, we get from above that for all large values of \,$k$, \,$\left\|\,s^{\,-\, 1}_{\,k}\,x \,-\, e,\, a_{\,2},\, \cdots,\, a_{\,n}\,\right\| \,<\, 1$\, and so by Corollary \ref{3.cor3.1}, \,$s^{\,-\, 1}_{\,k}\,x \,\in\, G$\, and therefore \,$x \,=\, s_{\,k}\,\left(\,s^{\,-\, 1}_{\,k}\,x\,\right)$\, belongs to \,$G$\, which contradicts the fact that \,$x \,\in\, S$.\,Thus, \,$\left\|\,s^{\,-\, 1}_{\,k},\, a_{\,2},\, \cdots,\, a_{\,n}\,\right\|$\, is not a bounded sequence.\,Clearly, we may assume that \,$\left\|\,s^{\,-\, 1}_{\,k},\, a_{\,2},\, \cdots,\, a_{\,n}\,\right\| \,\to\, \infty$\, as \,$k \,\to\, \infty$.\,Let \,$x_{\,k}$\, be defined as
\[x_{\,k} \,=\, \dfrac{s^{\,-\, 1}_{\,k}}{\left\|\,s^{\,-\, 1}_{\,k},\, a_{\,2},\, \cdots,\, a_{\,n}\,\right\|}\]
then \,$\left\|\,x_{\,k},\, a_{\,2},\, \cdots,\, a_{\,n}\,\right\| \,=\, 1$\, and
\begin{align*}
x\,x_{\,k}& \,=\, \dfrac{x\,s^{\,-\, 1}_{\,k}}{\left\|\,s^{\,-\, 1}_{\,k},\, a_{\,2},\, \cdots,\, a_{\,n}\,\right\|} \,=\, \dfrac{e \,+\, \left(\,x \,-\, s_{\,k}\,\right)\,s^{\,-\, 1}_{\,k}}{\left\|\,s^{\,-\, 1}_{\,k},\, a_{\,2},\, \cdots,\, a_{\,n}\,\right\|} \\
&=\,\dfrac{e}{\left\|\,s^{\,-\, 1}_{\,k},\, a_{\,2},\, \cdots,\, a_{\,n}\,\right\|} \,+\, \left(\,x \,-\, s_{\,k}\,\right)\,x_{\,k} \;\,\to\, \theta\; \;\text{as}\;\; \,k \,\to\, \infty\,, 
\end{align*}
because \,$\left\|\,s^{\,-\, 1}_{\,k},\, a_{\,2},\, \cdots,\, a_{\,n}\,\right\| \,\to\, \infty$, \,$s_{\,k} \,\to\, x$\, as \,$k \,\to\, \infty$, and \,$\left\|\,x_{\,k},\, a_{\,2},\, \cdots,\, a_{\,n}\,\right\| \,=\, 1$.\,Therefore, \,$x \,\in\, Z$\, and this proves the Theorem.      
\end{proof}

\subsection{Complex Homeomorphism in $n$-Banach Algebra}
\smallskip\hspace{.6 cm}

In this section, we introduce the idea of a complex homeomorphism in a \,$n$-Banach algebra and then deduce some of its properties.

\begin{definition}
Let \,$X$\, be a complex \,$n$-Banach algebra and \,$T$\, be a bounded \,$b$-linear functional defined on \,$X \,\times\, \left<\,b_{\,2}\,\right> \,\times\, \cdots \,\times\, \left<\,b_{\,n}\,\right>$.\,Then \,$T$\, is called a complex \,$b$-homeomorphism if
\begin{align*}
&T\,\left(\,x\,y,\, b_{\,2},\, \cdots,\, b_{\,n}\,\right) \,=\, T\,\left(\,x,\, b_{\,2},\, \cdots,\, b_{\,n}\,\right)\,T\,\left(\,y,\, b_{\,2},\, \cdots,\, b_{\,n}\,\right),
\end{align*}
for all \,$x,\, y \,\in\, X$.
\end{definition}

\begin{lemma}\label{3.lm3.31}
If \,$T$\, is complex \,$b$-homeomorphism defined on \,$X \,\times\, \left<\,b_{\,2}\,\right> \,\times\, \cdots \,\times\, \left<\,b_{\,n}\,\right>$, then \,$T\,\left(\,e,\, b_{\,2},\, \cdots,\, b_{\,n}\,\right) \,=\, 1$\, and if \,$x \,\in\, X$\, is invertible then \,$T\,\left(\,x,\, b_{\,2},\, \cdots,\, b_{\,n}\,\right) \,\neq\, 0$.
\end{lemma}

\begin{proof}
Since \,$T$\, is not identically zero, there exists \,$y$\, such that \,$T\,\left(\,y,\, b_{\,2},\, \cdots,\, b_{\,n}\,\right) \,\neq\, 0$.\,Now,
\begin{align*}
&T\,\left(\,y,\, b_{\,2},\, \cdots,\, b_{\,n}\,\right) \,=\, T\,\left(\,y\,e,\, b_{\,2},\, \cdots,\, b_{\,n}\,\right)=\,T\,\left(\,y,\, b_{\,2},\, \cdots,\, b_{\,n}\,\right)\,T\,\left(\,e,\, b_{\,2},\, \cdots,\, b_{\,n}\,\right)
\end{align*}
This implies that \,$T\,\left(\,e,\, b_{\,2},\, \cdots,\, b_{\,n}\,\right) \,=\, 1$.\,If \,$x$\, is invertible, then
\begin{align*}
T\,\left(\,x,\, b_{\,2},\, \cdots,\, b_{\,n}\,\right)\,T\,\left(\,x^{\,-\, 1},\, b_{\,2},\, \cdots,\, b_{\,n}\,\right)& \,=\, T\,\left(\,x\,x^{\,-\, 1},\, b_{\,2},\, \cdots,\, b_{\,n}\,\right)\\
&=\,T\,\left(\,e,\, b_{\,2},\, \cdots,\, b_{\,n}\,\right) \,=\, 1.
\end{align*} 
So, \,$T\,\left(\,x,\, b_{\,2},\, \cdots,\, b_{\,n}\,\right) \,\neq\, 0$.\,This proves the Lemma. 
\end{proof}

\begin{theorem}
Let \,$x \,\in\, X$\, and \,$\left\|\,x,\, a_{\,2},\, \cdots,\, a_{\,n}\,\right\| \,<\, 1$, for every \,$a_{\,2},\, \cdots,\, a_{\,n} \,\in\, X$.\,Then \,$\left|\,T\,\left(\,x,\, b_{\,2},\, \cdots,\, b_{\,n}\,\right)\,\right| \,<\, 1$, for every complex \,$b$-homeomorphism \,$T$\, defined on \,$X \,\times\, \left<\,b_{\,2}\,\right> \,\times\, \cdots \,\times\, \left<\,b_{\,n}\,\right>$.  
\end{theorem}

\begin{proof}
Let \,$\lambda$\, be a complex number such that \,$|\,\lambda\,| \,\geq\, 1$.\,Then \,$\left\|\,\dfrac{x}{\lambda},\, a_{\,2},\, \cdots,\, a_{\,n}\,\right\| \,<\, 1$\, and so by Theorem \ref{3.th3.11}, \,$e \,-\, \lambda^{\,-\, 1}\,x$\, is invertible.\,Therefore, by Lemma \ref{3.lm3.31}, \,$T\,\left(\,e \,-\, \lambda^{\,-\, 1}\,x,\, b_{\,2},\, \cdots,\, b_{\,n}\,\right) \,\neq\, 0$.\,But
\begin{align*}
T\,\left(\,e \,-\, \lambda^{\,-\, 1}\,x,\, b_{\,2},\, \cdots,\, b_{\,n}\,\right) &\,=\, T\,\left(\,e,\, b_{\,2},\, \cdots,\, b_{\,n}\,\right) \,-\, \lambda^{\,-\, 1}\,T\,\left(\,x,\, b_{\,2},\, \cdots,\, b_{\,n}\,\right)\\
&=\, 1 \,-\, \lambda^{\,-\, 1}\,T\,\left(\,x,\, b_{\,2},\, \cdots,\, b_{\,n}\,\right). 
\end{align*}
So, \,$T\,\left(\,x,\, b_{\,2},\, \cdots,\, b_{\,n}\,\right) \,\neq\, \lambda$\, which shows that \,$\left|\,T\,\left(\,x,\, b_{\,2},\, \cdots,\, b_{\,n}\,\right)\,\right| \,<\, 1$.\,This proves the Theorem. 
\end{proof}

In the following Theorem, we deal with the converse part of the Lemma \ref{3.lm3.31} i.\,e. we derive Gleason, Kahane, Zelazko type theorem with respect to complex \,$b$-homeomorphism.

\begin{theorem}
Let \,$T$\, be a bounded \,$b$-linear functional defined on \,$X \,\times\, \left<\,b_{\,2}\,\right> \,\times\, \cdots \,\times\, \left<\,b_{\,n}\,\right>$\, such that \,$T\,\left(\,e,\, b_{\,2},\, \cdots,\, b_{\,n}\,\right) \,=\, 1$\, and \,$T\,\left(\,x,\, b_{\,2},\, \cdots,\, b_{\,n}\,\right) \,\neq\, 0$, for every invertible \,$x \,\in\, X$.\,Then \,$T$\, is complex \,$b$-homeomorphism defined on \,$X \,\times\, \left<\,b_{\,2}\,\right> \,\times\, \cdots \,\times\, \left<\,b_{\,n}\,\right>$.   
\end{theorem}

\begin{proof}
Let \,$N$\, be the null space of \,$T$,\, i.\,e.,
\begin{align*}
N \,=\, \left\{\,x \,\in\, X \,:\, T\,\left(\,x,\, b_{\,2},\, \cdots,\, b_{\,n}\,\right) \,=\, 0\,\right\}. 
\end{align*}
Let \,$x \,\in\, X$\, and consider the element \,$x \,-\, \beta\,e$, where \,$\beta \,=\, T\,\left(\,x,\, b_{\,2},\, \cdots,\, b_{\,n}\,\right)$. Then
\begin{align*}
T\,\left(\,x \,-\, \beta\,e,\, b_{\,2},\, \cdots,\, b_{\,n}\,\right) &\,=\, T\,\left(\,x,\, b_{\,2},\, \cdots,\, b_{\,n}\,\right) \,-\, \beta\,T\,\left(\,e,\, b_{\,2},\, \cdots,\, b_{\,n}\,\right)\\
&=\,T\,\left(\,x,\, b_{\,2},\, \cdots,\, b_{\,n}\,\right) \,-\, \beta \,=\, 0
\end{align*}
and so \,$x \,-\, \beta\,e \,\in\, N$.\,Let \,$x \,-\, \beta\,e \,=\, a$.\,Then every element \,$x \,\in\, X$\, can be expressed as \,$x \,=\, a \,+\, T\,\left(\,x,\, b_{\,2},\, \cdots,\, b_{\,n}\,\right)\,e$, where \,$a \,\in\, N$.\,Similarly, if \,$y \,\in\, X$, let \,$y \,=\, b \,+\, T\,\left(\,y,\, b_{\,2},\, \cdots,\, b_{\,n}\,\right)\,e$, where \,$b \,\in\, N$.\,So,
\begin{align*}
x \,=\, a \,+\, T\,\left(\,x,\, b_{\,2},\, \cdots,\, b_{\,n}\,\right)\,e\,,\; \;y \,=\, b \,+\, T\,\left(\,y,\, b_{\,2},\, \cdots,\, b_{\,n}\,\right)\,e,
\end{align*} 
where \,$a,\, b \,\in\, N$\, and therefore
\begin{align*}
x\,y &\,=\, a\,b \,+\, T\,\left(\,y,\, b_{\,2},\, \cdots,\, b_{\,n}\,\right)\,a \,+\, T\,\left(\,x,\, b_{\,2},\, \cdots,\, b_{\,n}\,\right)\,b \,+\\
&\hspace{1cm} \,+\, T\,\left(\,x,\, b_{\,2},\, \cdots,\, b_{\,n}\,\right)\,T\,\left(\,y,\, b_{\,2},\, \cdots,\, b_{\,n}\,\right)\,e.  
\end{align*}
Thus,
\begin{align}
T\,\left(\,x\,y,\, b_{\,2},\, \cdots,\, b_{\,n}\,\right)& \,=\, T\,\left(\,x,\, b_{\,2},\, \cdots,\, b_{\,n}\,\right)\,T\,\left(\,y,\, b_{\,2},\, \cdots,\, b_{\,n}\,\right) \,+\nonumber\\
&\hspace{1cm}+\,T\,\left(\,a\,b,\, b_{\,2},\, \cdots,\, b_{\,n}\,\right).\label{3.eq3.38}
\end{align}
Therefore, the theorem will be proved if \,$T\,\left(\,a\,b,\, b_{\,2},\, \cdots,\, b_{\,n}\,\right) \,=\, 0$\, i.\,e., if 
\begin{align}
a\,b \,\in\, N\; \;\text{whenever}\; \;a \,\in\, N\; \;\text{and}\; \;b \,\in\, N.\label{3.eq3.39}
\end{align}
Suppose that a special case of (\ref{3.eq3.39}) is true i.\,e.,
\begin{align}
a^{\,2} \,\in\, N\; \;\text{if}\; \; a \,\in\, N.\label{3.eq3.40}
\end{align}
Assuming (\ref{3.eq3.40}), we can deduce (\ref{3.eq3.39}) as follows.\,In (\ref{3.eq3.38}), we assume \,$x \,=\, y$\, and so \,$a \,=\, b$.\,Then, because of (\ref{3.eq3.40}), we get
\begin{align}
T\,\left(\,x^{\,2},\, b_{\,2},\, \cdots,\, b_{\,n}\,\right) &\,=\, T\,\left(\,a^{\,2},\, b_{\,2},\, \cdots,\, b_{\,n}\,\right) \,+\, \left[\,T\,\left(\,x,\, b_{\,2},\, \cdots,\, b_{\,n}\,\right)\,\right]^{\,2} \nonumber\\
&=\,\left[\,T\,\left(\,x,\, b_{\,2},\, \cdots,\, b_{\,n}\,\right)\,\right]^{\,2}\,,\; x \,\in\, X.\label{3.eq3.41}
\end{align} 
In (\ref{3.eq3.41}), we replace \,$x$\, by \,$x \,+\, y$\, and get
\begin{align*}
&T\,\left(\,x\,y \,+\, y\,x,\, b_{\,2},\, \cdots,\, b_{\,n}\,\right) \,=\, 2\,T\,\left(\,x,\, b_{\,2},\, \cdots,\, b_{\,n}\,\right)\,T\,\left(\,y,\, b_{\,2},\, \cdots,\, b_{\,n}\,\right),\, x,\, y \,\in\, X.
\end{align*}
Therefore, 
\begin{align}
x\,y \,+\, y\,x \,\in\, N\;\; \;\;\text{if}\; \,x \,\in\, N\,,\; \;y \,\in\, X.\label{3.eq3.42}
\end{align}
By a simple calculation we obtain
\begin{align}
\left(\,x\,y \,-\, y\,x\,\right)^{\,2} \,+\, \left(\,x\,y \,+\, y\,x\,\right)^{\,2} \,=\, 2\,\left[\,x\,(\,y\,x\,y\,) \,+\, (\,y\,x\,y\,)\,x\,\right].\label{3.eq3.43} 
\end{align}
If \,$x \,\in\, N$, then the right side of (\ref{3.eq3.43}) belongs to \,$N$\, because of (\ref{3.eq3.42}).\,Also, for the same reason \,$x\,y \,+\, y\,x \,\in\, N$\, and so by an application of (\ref{3.eq3.41}), \,$ \left(\,x\,y \,+\, y\,x\,\right)^{\,2} \,\in\, N$. Therefore, \,$ \left(\,x\,y \,-\, y\,x\,\right)^{\,2} \,\in\, N$\, and another application of (\ref{3.eq3.41}) gives that
\begin{align}
x\,y \,-\, y\,x \,\in\, N\;\; \;\;\text{if}\; \,x \,\in\, N\,,\; \;y \,\in\, X.\label{3.eq3.44}
\end{align}
Adding (\ref{3.eq3.42}) and (\ref{3.eq3.44}), we get (\ref{3.eq3.39}) and so the theorem is proved provided we establish  (\ref{3.eq3.40}) which we do now.

Since by hypothesis \,$T\,\left(\,x,\, b_{\,2},\, \cdots,\, b_{\,n}\,\right) \,\neq\, 0$, for every invertible \,$x \,\in\, X$, \,$N$\, contains no invertible element.\,So, if \,$\alpha \,\in\, N$, then by Corollary \ref{3.cor3.1}, 
\[\left\|\,e \,-\, \alpha,\, b_{\,2},\, \cdots,\, b_{\,n}\,\right\| \,\geq\, 1.\]
Clearly, \,$\dfrac{\alpha}{\lambda} \,\in\, N$, for any complex number \,$\lambda$\, and so \,$\left\|\,e \,-\, \dfrac{\alpha}{\lambda},\, b_{\,2},\, \cdots,\, b_{\,n}\,\right\| \,\geq\, 1$.\,Therefore,
\begin{align}
\left\|\,\lambda\,e \,-\, \alpha,\, b_{\,2},\, \cdots,\, b_{\,n}\,\right\|& \,=\, |\,\lambda\,|\,\left\|\,e \,-\, \dfrac{\alpha}{\lambda},\, b_{\,2},\, \cdots,\, b_{\,n}\,\right\| \nonumber\\
&\geq\, |\,\lambda\,| \,=\, \left|\,T\,\left(\,\lambda\,e \,-\, \alpha,\, b_{\,2},\, \cdots,\, b_{\,n}\,\right)\,\right|.\label{3.eq3.45}
\end{align} 
As in the first part of the proof, we can show that every element \,$x \,\in\, X$\, can be expressed as \,$x \,=\, \lambda\,e\, \,-\, \alpha$, where \,$\alpha \,\in\, N$\, and \,$\lambda$\, is a complex number.\,So, by (\ref{3.eq3.45}), we get that
\[\left|\,T\,\left(\,x,\, b_{\,2},\, \cdots,\, b_{\,n}\,\right)\,\right| \,\leq\, \left\|\,x,\, b_{\,2},\, \cdots,\, b_{\,n}\,\right\|,\]
for all \,$x \,\in\, X$.\,Thus, \,$T$\, is a bounded \,$b$-linear functional defined on \,$X \,\times\, \left<\,b_{\,2}\,\right> \,\times\, \cdots \,\times\, \left<\,b_{\,n}\,\right>$\, with norm \,$1$.

Let \,$a \,\in\, N$\, and without loss of generality, we may assume that \,$\left\|\,a,\, b_{\,2},\, \cdots,\, b_{\,n}\,\right\| \,=\, 1$.\,Now, we define
\begin{align*}
f\,(\,\lambda\,) \,=\, \sum\limits_{j \,=\, 0}^{\,\infty}\,\dfrac{T\,\left(\,a^{\,j},\, b_{\,2},\, \cdots,\, b_{\,n}\,\right)}{j\,!}\,\lambda^{\,j},
\end{align*}
where \,$\lambda$\, is a complex number.\,Now,
\begin{align*}
\left|\,T\,\left(\,a^{\,j},\, b_{\,2},\, \cdots,\, b_{\,n}\,\right)\,\right| \,\leq\, \left\|\,a^{\,j},\, b_{\,2},\, \cdots,\, b_{\,n}\,\right\| \,\leq\, \left\|\,a,\, b_{\,2},\, \cdots,\, b_{\,n}\,\right\|^{\,j} \,=\, 1, 
\end{align*} 
for all \,$j$\, and so \,$f$\, is a entire function and further 
\begin{align*}
\left|\,f\,(\,\lambda\,)\,\right| \,\leq\, \sum\limits_{j \,=\, 0}^{\,\infty}\,\dfrac{\left|\,T\,\left(\,a^{\,j},\, b_{\,2},\, \cdots,\, b_{\,n}\,\right)\,\right|}{j\,!}\,|\,\lambda\,|^{\,j} \,\leq\, \sum\limits_{j \,=\, 0}^{\,\infty}\,\dfrac{|\,\lambda\,|^{\,j}}{j\,!} \,=\, e^{\,|\,\lambda\,|}. 
\end{align*}
Also, by hypothesis \,$f\,(\,0\,) \,=\, T\,\left(\,e,\, b_{\,2},\, \cdots,\, b_{\,n}\,\right) \,=\, 1$\, and since \,$a \,\in\, N$, we get \,$f^{\,\prime}\,(\,0\,) \,=\, T\,\left(\,a,\, b_{\,2},\, \cdots,\, b_{\,n}\,\right) \,=\, 0$.

Thus, if we can prove that \,$\left|\,f\,(\,\lambda\,)\,\right| \,>\, 0$\, for every complex number \,$\lambda$, then Lemma \ref{3.lm3.31} gives that \,$f\,(\,\lambda\,) \,=\, 1$, for all \,$\lambda$\, i.\,e., \,$f^{\,\prime\,\prime}\,(\,0\,) \,=\, 0$.\,But \,$f^{\,\prime\,\prime}\,(\,0\,) \,=\, T\,\left(\,a^{\,2},\, b_{\,2},\, \cdots,\, b_{\,n}\,\right)$\, and so \,$T\,\left(\,a^{\,2},\, b_{\,2},\, \cdots,\, b_{\,n}\,\right) \,=\, 0$\, i.\,e., \,$a^{\,2} \,\in\, N$\, and this is (\ref{3.eq3.40}).

Consider the series \,$L\,(\,\lambda\,) \,=\, \sum\limits_{j \,=\, 0}^{\,\infty}\,\dfrac{\lambda^{\,j}}{j\,!}\,a^{\,j}$, where \,$\lambda$\, is a complex number and this is converges in the \,$n$-norm of \,$X$.\,Since \,$T$\, is \,$b$-bounded, we obtain
\begin{align}
T\,\left(\,L\,(\,\lambda\,),\, b_{\,2},\, \cdots,\, b_{\,n}\,\right)& \,=\, T\,\left(\,\sum\limits_{j \,=\, 0}^{\,\infty}\,\dfrac{\lambda^{\,j}}{j\,!}\,a^{\,j},\, b_{\,2},\, \cdots,\, b_{\,n}\,\right)\nonumber\\ 
&=\,\sum\limits_{j \,=\, 0}^{\,\infty}\,\dfrac{T\,\left(\,a^{\,j},\, b_{\,2},\, \cdots,\, b_{\,n}\,\right)}{j\,!}\,\lambda^{\,j} \,=\, f\,(\,\lambda\,).\label{3.eq3.46}
\end{align}
From the expression for \,$L\,(\,\lambda\,)$, we can verify, as in the classical case, the following functional equation \,$L\,(\,\lambda \,+\, \mu\,) \,=\, L\,(\,\lambda\,)\,L\,(\,\mu\,)$.\,We have in particular \,$L\,(\,\lambda\,)\,L\,(\, \,-\, \lambda\,) \,=\, L\,(\,0\,) \,=\, e$\, which shows that \,$L\,(\,\lambda\,)$\, is an invertible element in \,$X$.\,By hypothesis then \,$T\,\left(\,L\,(\,\lambda\,),\, b_{\,2},\, \cdots,\, b_{\,n}\,\right) \,\neq\, 0$, i.\,e., from (\ref{3.eq3.46}), we get \,$f\,(\,\lambda\,) \,\neq\, 0$.\,This completes the proof.      
\end{proof}

\end{document}